\documentclass[a4paper,reqno, oneside]{amsart}
\usepackage[english]{babel}
\usepackage{amsmath, amssymb, amsthm, amscd}
\usepackage{enumerate}
\usepackage{palatino}
\usepackage{mathpazo}
\usepackage{paralist}
\usepackage[a4paper]{geometry}

\usepackage{tikz} 
\usepackage{tikz-cd}

 \usepackage{amsaddr}

\DeclareMathOperator{\id}{id}

\newcommand{\NN}{\mathbb{N}}
\newcommand{\ZZ}{\mathbb{Z}}

\newcommand{\CC}{\mathbb{C}}
\newcommand{\HH}{\mathbb{H}}
\newcommand{\KK}{\mathbb{K}}
\newcommand{\QQ}{\mathbb{Q}}

\newcommand{\cpr}{\vee}

\makeatletter
\newcommand{\ostar}{\mathbin{\mathpalette\make@circled\ast}}
\newcommand{\make@circled}[2]{%
  \ooalign{$\m@th#1\smallbigcirc{#1}$\cr\hidewidth$\m@th#1#2$\hidewidth\cr}%
}
\newcommand{\smallbigcirc}[1]{%
  \vcenter{\hbox{\scalebox{0.77778}{$\m@th#1\bigcirc$}}}%
}
\makeatother

\renewcommand{\setminus}{\smallsetminus}

\theoremstyle{plain}
\newtheorem{theorem}{Theorem}[section]
\newtheorem{prop}[theorem]{Proposition}
\newtheorem{lemma}[theorem]{Lemma}
\newtheorem{cor}[theorem]{Corollary}  
\newtheorem*{claim}{Claim}
\newtheorem*{thm}{Theorem}

\theoremstyle{definition}
\newtheorem{definition}[theorem]{Definition}

\theoremstyle{remark}
\newtheorem{remark}[theorem]{Remark}
\newtheorem{example}[theorem]{Example}

\numberwithin{equation}{section}

\begin{document}
\setlength{\parindent}{0.cm}

\title{The string topology coproduct on complex and quaternionic projective space}

\author{Maximilian Stegemeyer}
\address{Mathematisches Institut, Universit\"at Freiburg, Ernst-Zermelo-Straße 1, 79104 Freiburg, Germany}
\email{maximilian.stegemeyer@math.uni-freiburg.de}
\date{\today}
\maketitle

\begin{abstract}
On the free loop space of compact symmetric spaces Ziller introduced explicit cycles generating the homology of the free loop space.
We use these explicit cycles to compute the string topology coproduct on complex and quaternionic projective space.
The behavior of the Goresky-Hingston product for these spaces then follows directly.
\end{abstract}

\setcounter{tocdepth}{1} 
\tableofcontents

\section{Introduction}

The central idea in Morse theory is to study the interaction between the critical sets of a function on a differentiable manifold and the topology of this manifold.
While it is usually easy to understand the local homology around a critical level, it is a hard question to determine if and how all of the homology of the manifold can be understood by the individual homologies of the critical sets.
In \cite{ziller:1977} Ziller defines cycles on the free loop space of a compact globally symmetric space which can be used to show that the relative cycles from level homology can be completed in the free loop space.
This idea goes back to Bott's $K$-cycles, see \cite{bott:1956} and \cite{bott:1958b} as well as Bott's and Samelson's work in \cite{bott:1958a}.
Hingston and Oancea use explicit cycles in the path space of complex projective space to compute a Pontryagin-Chas-Sullivan type product in \cite{hingston:2013oancea}.
There they use the name \textit{completing manifold} for this construction of completing relative cycles.
We shall use this terminology here as well.
The goal of this article is to use Ziller's completing manifolds to compute the string topology coproduct for complex and quaternionic projective space.

The string topology coproduct was introduced by Goresky and Hingston in \cite{goresky:2009}.
It was furthermore studied by Hingston and Wahl in \cite{hingston:2017} where the authors give a definition which is equivalent to the one in \cite{goresky:2009} and which we shall use in this article.
If $M$ is a closed oriented Riemannian manifold of dimension $N$, the string topology coproduct is a map
$$   \cpr \colon \mathrm{H}_i(\Lambda M,M)\to \mathrm{H}_{i+1-N}(\Lambda M\times \Lambda M,\Lambda M\times M\cup M\times \Lambda M)  $$
where $\Lambda M$ is the free loop space of $M$ and $M$ is considered as a subspace of $\Lambda M$ via the identification of a point with the trivial loop at this point.
The string topology coproduct has been computed for odd-dimensional spheres in \cite{hingston:2017}.
Furthermore, there are partial computations of the string topology coproduct on Lens spaces, see \cite{naef:2021} and \cite{naef:2022riveraWahl}.
In particular in \cite{naef:2022riveraWahl} Naef, Rivera and Wahl show that the string topology coproduct is not a homotopy invariant in general.
However, as Hingston and Wahl show in \cite{hingston:2019} if one only considers homotopy invariances with certain additional conditions then the string topology coproduct is invariant under these maps.

The author of this article used Bott's $K$-cycles - which can be understood as completing manifolds - to show that the string topology coproduct is trivial for compact simply connected Lie groups of rank $r\geq 2$, see \cite{stegemeyer:2021}.
In this present article we use completing manifolds to compute the string topology coproduct on $M =\CC P^n$ and $M =\HH P^n$.
We are going to explicitly describe a family of closed manifolds $\Gamma_k$ and embeddings $f_k\colon \Gamma_k\to \Lambda M$, $k\in \NN$ such that there is a family of homology classes
$$  A_k^i = (f_k)_*[\alpha_k^i] \in \mathrm{H}_{\bullet}(\Lambda M,M;\QQ) \quad \text{and}\quad B_k^i =(f_k)_*[\beta_k^i]\in \mathrm{H}_{\bullet}(\Lambda M,M;\QQ)  ,$$ $ k\in \NN,i\in\{0,\ldots,n-1\}        $
which generate all of the homology of the pair $(\Lambda M,M)$ with rational coefficients.
Here $[\alpha_k^i]$ and $[\beta_k^i]$ are classes in the homology of $\Gamma_k$.
We will then show that the string topology coproduct behaves as follows.
\begin{thm}[Theorem \ref{theorem_coproduct}]
Let $\mathbb{K}$ be either $\mathbb{C}$ or $\mathbb{H}$.
The string topology coproduct on $M = \KK P^n$ satisfies
$$  \cpr A_k^i = \sum_{m=1}^{k-1} \sum_{j=0}^i A_m^j\times A_{k-m}^{i-j}           $$
and
$$    \cpr B_k^i = \sum_{m=1}^{k-1} \sum_{j=0}^i ( B_m^j \times A_{k-m}^{i-j} - A_m^j \times B_{k-m}^{i-j}  ) .       $$
\end{thm}
We use this result to compute the Goresky-Hingston product on the manifolds $\CC P^n$ and $\HH P^n$, see Theorem \ref{theorem_gh_product}.

This article is organized as follows.
In Section \ref{sec_completing_manifolds} we introduce the notion of completing manifolds and discuss their relevance in Morse theory.
The string topology coproduct is defined in Section \ref{sec_string_topology_cpr}.
We study the critical manifolds of the length functional in the free loop space and in the figure eight space in Section \ref{sec_k_cycles}.
The completing manifolds are introduced in Section \ref{sec_construction_completing} and in Section \ref{sec_cohomology_completing} we discuss their cohomology ring.
The computation of the string topology coproduct is then carried out in Section \ref{sec_computation}.
Finally, in Section \ref{sec_cohomology} we use the results of the previous sections to compute the Goresky-Hingston product on $\CC P^n$ and $\HH P^n$. 

The proof of a central Lemma of Section \ref{sec_computation} is to be found in Appendix \ref{appendix_thom_class} and in Appendix \ref{appendix_cap} we discuss a relative version of the standard cap product which is used in the definition of the string topology coproduct.

\subsection*{Acknowledgements}

The author thanks the anonymous referee for their careful and thoughtful reading of our manuscript. Their suggestions highly improved the exposition of this article.

\section{Completing manifolds and Morse theory} \label{sec_completing_manifolds}

We start by introducing the notion of a completing manifold following the expositions in \cite{hingston:2013oancea} and \cite{oancea:2015}.

Let $X$ be a Hilbert manifold and let $f \colon \mathbb{X}\to \mathbb{R}$ be a smooth function on $X$ satisfying the Palais-Smale condition (C).
Let $a$ be a critical value of $f$ and assume that the set of critical points $B$ at level $a$ is a non-degenerate finite-dimensional critical submanifold of finite index $k$ with orientable negative bundle.
Then the behavior of the level homology $\mathrm{H}_{\bullet}(X^{\leq a},X^{<a})$ is well known.
It holds that 
$$   \mathrm{H}_{\bullet}(X^{\leq a},X^{<a}) \cong \mathrm{H}_{\bullet -k}(B)        $$
where coefficients can be taken in an arbitrary commutative ring $R$.
In applications the homology of these critical submanifolds may be much easier to understand than the homology of $X$.
Therefore, one would like to find conditions which imply that all of the homology of $X$ is built up by these level homologies.
\begin{definition}[\cite{oancea:2015}, Definition 6.1] \label{def_completing_mfld}
Let $X$ be a Hilbert manifold and let $f$ be a smooth real-valued function on $X$ satisfying condition (C).
Let $a$ be a critical value of $f$ and assume that $B$ is a non-degenerate connected critical submanifold at level $a$ of index $k$ and of dimension $l = \mathrm{dim}(B)$.
Assume that $k$ and $l$ are both finite.
A \textit{completing manifold} for $B$ is a closed, orientable manifold $\Gamma$ of dimension $k + l$ with an embedding $\varphi: \Gamma \to X^{\leq a}$ such that the following holds.
There is an $l$-dimensional submanifold $L$ such that $\varphi|_L$ maps $L$ homeomorphically onto $B$ and there is a retraction map $p:\Gamma\to L$.
Furthermore, the embedding $\varphi$ induces a map of pairs
$$  \varphi \colon (\Gamma,\Gamma\setminus L) \to (X^{\leq a}, X^{<a}) .    $$
\end{definition}

\begin{remark}
\begin{enumerate}
    \item 
This definition of a completing manifold is actually the one of a \textit{strong} completing manifold in \cite{oancea:2015}.
Since all cases that we consider in this article satisfy the assumption of this stronger version we limit our attention to this situation.
\item Note that the above definition can be used for cases where the critical set at level $a$ consists of several connected critical submanifolds.
We can then set up a completing manifold for each connected component.
We will see this in the case of the figure-eight space in Section \ref{sec_k_cycles}.
\end{enumerate}
\end{remark}

Recall that if $f\colon M\to N$ is a map between oriented manifolds then the Gysin map
$$f_!\colon \mathrm{H}_j(N)\to \mathrm{H}_{j+\mathrm{dim}(M)-\mathrm{dim}(N)}(M)$$ is given by
$$f_! \colon \mathrm{H }_j(N) \xrightarrow[]{(PD_{N})^{-1}}  \mathrm{H}^{\mathrm{dim}(N) - j}(N) \xrightarrow[]{f^*} \mathrm{H}^{\mathrm{dim}(N)-j}(M) \xrightarrow[]{PD_{M}} \mathrm{H}_{\mathrm{dim}(M) - (\mathrm{dim}(N)-j)}(M) . $$
Here $PD_B$ stands Poincaré duality on the manifold $B$. 
The Gysin map $p_!\colon \mathrm{H}_i(L)\to \mathrm{H}_{i+k}(\Gamma)$ is clearly a right inverse to the Gysin map $s_!\colon \mathrm{H}_i(\Gamma)\to \mathrm{H}_{i-k}(L)$ where $s\colon L\hookrightarrow \Gamma$ is the embedding of $L$ into $\Gamma$ given by the data of the completing manifold.
Up to sign, the Gysin map $s_!$ is equal to the composition
$$    \mathrm{H}_i(\Gamma) \to \mathrm{H}_i(\Gamma,\Gamma\setminus L) \xrightarrow[]{\cong}  \mathrm{H}_{i-k}(L)   $$
where the first map is induced by the inclusion of pairs and the second is the Thom isomorphism, see \cite[Theorem VI.11.3]{bredon:2013}.
This shows that the map $\mathrm{H}_{\bullet}(\Gamma)\to \mathrm{H}_{\bullet}(\Gamma,\Gamma\setminus L)$ is surjective.
See also \cite[Remark 7]{hingston:2013oancea}.
In particular this observation leads to the following result.

\begin{prop}[\cite{oancea:2015}, Lemma 6.2] \label{prop_completing_mfld}
Let $X$ be a Hilbert manifold, $f$ a smooth real-valued function on $X$ satisfying condition (C) and $a$ be a critical value of $f$.
Assume that the set of critical points at level $a$ is a non-degenerate critical submanifold $B$ of index $k$.
If there is a completing manifold for $B$ then
$$   \mathrm{H}_{\bullet}(X^{\leq a}) \cong \mathrm{H}_{\bullet}(X^{<a}) \oplus \mathrm{H}_{\bullet }(X^{\leq a},X^{<a}) \cong \mathrm{H}_{\bullet}(X^{<a}) \oplus \mathrm{H}_{\bullet -k}(B)  .         $$
\end{prop}

If the homology of the sublevel set $X^{\leq a}$ is isomorphic to the direct sum
\begin{equation} \label{eq_perfect_morse}
        \mathrm{H}_{\bullet}(X^{\leq a}) \cong \mathrm{H}_{\bullet}(X^{<a}) \oplus \mathrm{H}_{\bullet}(X^{\leq a},X^{<a})       
\end{equation}
for all critical values $a$, we say that the function $f$ is a \textit{perfect} Morse-Bott function.
This property clearly holds if all the connecting morphisms in the long exact sequence of the pair $(X^{\leq a},X^{<a})$ vanish.
If every critical submanifold has a completing manifold, it follows that the function $f$ is perfect.
Using completing manifolds Ziller shows in \cite{ziller:1977} that the energy function on the free loop space of a compact symmetric space is a perfect Morse-Bott function.
Note that he uses $\ZZ_2$-coefficients in general, since there are issues with orientability for some spaces.
We will describe these completing manifolds for $M = \CC P^n$ and $M = \HH P^n$ in detail in Section \ref{sec_k_cycles}.

There is also an obvious generalization of the above Proposition if we are in the situation of the critical set decomposing into several connected components and each one admitting a completing manifold.

\section{The string topology coproduct} \label{sec_string_topology_cpr}

In this section we introduce the string topology coproduct.
We closely follow the definition of the coproduct given in \cite{hingston:2017}.
Let $M$ be an oriented closed $N$-dimensional Riemannian manifold.
We denote the unit interval by $I = [0,1]$.
Let $$  PM = \big\{ \gamma : I\to M\,|\, \gamma \,\,\text{absolutely continuous}, \,\, \int_0^1 |\Dot{\gamma}(t)|^2 \,\mathrm{d}t < \infty \big\}   $$
be the set of absolutely continuous curves in $M$ such that their derivative is square integrable.
See \cite[Definition 2.3.1]{klingenberg:1995} for the definition of absolutely continuous curves in a manifold.
We define the free loop space of $M$ to be
$$  \Lambda M = \{\gamma \in PM \,|\, \gamma(0 ) = \gamma(1)\}   $$
and this is in fact a submanifold of $PM$.
The manifold $M$ itself can be embedded into $\Lambda M$ via the trivial loops, see \cite[Proposition 1.4.6]{klingenberg:78}.
On the path space $PM$ we consider the length functional
\begin{equation} \label{eq_length_functional}
      \mathcal{L}: PM \to [0,\infty),\qquad  \mathcal{L}(\gamma) = \sqrt{   \int_0^1 |\Dot{\gamma}(t) |^2 \mathrm{d}t   }  
\end{equation}
which is a continuous function on $PM$, see \cite[Theorem 2.3.20]{klingenberg:1995}.
Moreover, it is smooth on $PM\setminus M$.
If we restrict $\mathcal{L}$ to the free loop space $\Lambda M$ it turns out that the non-trivial critical points of $\mathcal{L}$ are precisely the closed geodesics in $M$.

We now fix a commutative ring $R$ and consider homology and cohomology with coefficients in $R$.
Fix an $\epsilon>0$ smaller than the injectivity radius of $M$.
Then the diagonal $\Delta M\subseteq M\times M$ has a tubular neighborhood given by
$$   U_M = \{ (p,q)\in M\times M \,|\, \mathrm{d}(p,q) < \epsilon \}  .  $$
Here, $\mathrm{d}$ is the distance function on $M$ induced by the Riemannian metric.
For an $\epsilon_0> 0$ such that $\epsilon_0<\epsilon$ we set
$$  U_{M,\geq \epsilon_0} = \{(p,q)\in U_M\,|\, \mathrm{d}(p,q)\geq \epsilon_0\} .    $$
As Hingston and Wahl argue, see \cite[Section 1.3]{hingston:2017}, the Thom class in $\mathrm{H}^{N}(TM,TM\setminus M)$ induces a Thom class $\tau_M\in \mathrm{H}^{N}(U_M,U_{M,\geq\epsilon_0})$.
On the free loop space we consider the space
$$ F_{\Lambda} = \{(\gamma,s)\in \Lambda M\times I\,|\,\gamma(s) = \gamma(0)\} .  $$
We set
\begin{eqnarray*}
  U_{\Lambda} &=& \{ (\gamma,s) \in \Lambda M\times I\,|\, \mathrm{d}(\gamma(0),\gamma(s))< \epsilon\}  \qquad \text{and} \\
  U_{\Lambda,\geq\epsilon_0} & = &
  \{ (\gamma,s)\in U_{\Lambda}\,|\, \mathrm{d}(\gamma(0),\gamma(s))\geq \epsilon_0\} .
\end{eqnarray*}
The set $U_{\Lambda}$ is an open neighborhood of $F_{\Lambda}$.
Define the evaluation map $\mathrm{ev}_{\Lambda}: \Lambda\times I\to M\times M$ by
$$ \mathrm{ev}_{\Lambda}(\gamma,s) = (\gamma(0),\gamma(s)) . $$
This yields a map of pairs
$$   \mathrm{ev}_{\Lambda}\colon (U_{\Lambda},U_{\Lambda,\geq\epsilon_0})   \to (U_M,U_{M,\geq\epsilon_0}) . $$
We define the class
$$  \tau_{\Lambda} = \mathrm{ev}_{\Lambda}^* \tau_M \in \mathrm{H}^N(U_{\Lambda},U_{\Lambda,\geq\epsilon_0}) .    $$
Furthermore, there is a retraction map $$   R_{GH}\colon U_{\Lambda}\to F_{\Lambda}  . $$
We refer to \cite[Section 1.5]{hingston:2017} for its precise definition.
Finally, consider the cutting map
$$  \mathrm{cut}\colon F_{\Lambda}\to \Lambda M\times \Lambda M      $$
which maps a point $(\gamma,s)\in\Lambda\times I$ with $\gamma(0)=\gamma(s)$ to the pair of loops
$   (\gamma|_{[0,s]},\gamma|_{[s,1]})  $ and reparametrizes both loops such that they are again defined on the unit interval $I$.
Note that the cutting map actually factors through maps
$$  F_{\Lambda} \xrightarrow[]{\widetilde{\mathrm{cut}}} \Lambda M\times_M \Lambda M \hookrightarrow \Lambda M\times \Lambda M     $$
where $\Lambda M\times_M \Lambda M$ is the figure-eight space
$$ \Lambda M\times_M \Lambda M = \{ (\gamma,\sigma)\in \Lambda M\times \Lambda M\,|\, \gamma(0) = \sigma(0) \}  . $$

With this preparation we can now define the string topology coproduct.
Let $[I]$ be the positively oriented generator of $\mathrm{H}_1(I,\partial I)$ with respect to the standard orientation of the unit interval.
In order to shorten notation we shall also write $\Lambda$ for the free loop space $\Lambda M$.

\begin{definition}
The \textit{string topology coproduct} is defined as the map
\begin{eqnarray*}
 \cpr: \mathrm{H}_{\bullet}(\Lambda,M) &\xrightarrow{\times [I]}& \mathrm{H}_{\bullet+1}(\Lambda\times I, \Lambda\times\partial I\cup M\times I)
     \\ &\xrightarrow{\tau_{\Lambda}\cap }& \mathrm{H}_{\bullet+1-N}(U_{\Lambda}, \Lambda\times\partial I\cup M\times I) \\
    &\xrightarrow{(\mathrm{R}_{GH})_*}&
    \mathrm{H}_{\bullet+1-N}(F_{\Lambda}, \Lambda\times\partial I\cup M\times I) \\
    &\xrightarrow{(\mathrm{cut})_*}& \mathrm{H}_{\bullet +1-N}(\Lambda\times\Lambda, \Lambda\times M\cup M\times \Lambda). 
    \end{eqnarray*}
\end{definition}
\begin{remark} Let $M$ be a closed oriented manifold.
\begin{enumerate}
    \item Note that the cap product with the class $\tau_{\Lambda}$ is a particular relative cap product.
    This relative cap product is defined in Appendix \ref{appendix_cap} where we also study some basic properties. See also \cite[Appendix A]{hingston:2017}.
    \item Hingston and Wahl define an \textit{algebraic loop coproduct}, see \cite[Definition 1.6]{hingston:2017}, which is a sign-corrected version of the string topology coproduct.
    Since we will later only consider even-dimensional manifolds, this sign correction does not matter.
\end{enumerate}
\end{remark}
If we use field coefficients, then the string topology coproduct induces a dual product in cohomology which is known as the Goresky-Hingston product.
\begin{definition}
Let $\mathbb{F}$ be a field and assume that the homology of $\Lambda M$ is of finite type.
Let $\alpha\in \mathrm{H}^i(\Lambda ,M;\mathbb{F})$ and $\beta\in \mathrm{H}^j(\Lambda ,M ;\mathbb{F})$ be relative cohomology classes, then the \textit{Goresky-Hingston product} $\alpha\ostar \beta$ is defined to be the unique cohomology class in $\mathrm{H}^{i+j+N-1}(\Lambda,M ;\mathbb{F})$ such that
$$   \langle \alpha \ostar \beta , X\rangle = \langle \alpha \times \beta , \cpr X\rangle \qquad \text{for all} \,\,\, X\in\mathrm{H}_{\bullet}(\Lambda M,M;\mathbb{F}).    $$
\end{definition}
\begin{remark} Let $M$ be a closed oriented manifold.
\begin{enumerate}
    \item The Goresky-Hingston product can also be defined intrinsically, see \cite{goresky:2009}. However, in this article we shall study properties of this product only via the duality with the string topology coproduct.
    \item As for the string topology coproduct, Hingston and Wahl define a sign-corrected version of the Goresky-Hingston product in \cite{hingston:2017}.
    For even-dimensional manifolds, the above product and its sign-corrected version agree.
    Hence, in this article the distinction will not matter.
\end{enumerate}
\end{remark}

\begin{remark} \label{rem_sum}
Since the next four sections will deal with the technical details of the computation of the coproduct we want to sum up the strategy for computing the string topology coproduct on $ M = \CC P^n$ or $M = \HH P^n$ at this point.
\begin{itemize}
    \item In Section \ref{sec_k_cycles} we will study the critical manifolds $\Sigma_k$, $k\in \NN$ in $\Lambda M$ of the length functional $\mathcal{L}\colon \Lambda M\to \mathbb{R}$.
    \item In Section \ref{sec_construction_completing} we construct the completing manifolds $\Gamma_k$.
    \item We shall see that the manifold $\Gamma_k$, $k\in \NN$ can also serve as a completing manifold for critical submanifolds in the figure-eight space $\Lambda M\times_M\Lambda M$.
    \item We determine the cohomology ring of $\Gamma_k$ in Section \ref{sec_cohomology_completing}. We can then explicitly compute the Gysin map and obtain a set of generators for the homology $\mathrm{H}_{\bullet}(\Lambda M,M)$.
    \item In Section \ref{sec_computation} we will then replicate all the steps in the definition of the coproduct on the manifold $\Gamma_k$.
    \item We will pull back the class $\tau_{\Lambda}$ to a class which can be described in terms of the cohomology of $\Gamma_k$ and compute the cap product with this class.
    \item Then one sees that under the cutting map $\widetilde{\mathrm{cut}}\colon F_{\Lambda}\to \Lambda M \times_M\Lambda M$ we get homology classes which we can identify with classes coming from the manifold $\Gamma_k$ seen as a completing manifold in the figure-eight space $\Lambda M\times_M \Lambda M$.
\end{itemize}
\end{remark}

\section{Critical manifolds in the free loop space of projective spaces} \label{sec_k_cycles}

In this section we describe the completing manifolds on the loop space of $\mathbb{C} P^n$ and $\HH P^n$ in detail.
We will first study the critical manifolds of $\Lambda M$ and $\Lambda M\times_M \Lambda M$ with respect to the length functional $\mathcal{L}$.
Then we define the completing manifolds $\Gamma_k$, $k\in \NN$ and show that $\Gamma_k$ can serve as a completing manifold both in $\Lambda M$ as well as in $\Lambda M\times_M \Lambda M$.
Finally, we describe the cohomology ring of $\Gamma_k$ in detail and give a set of explicit generators of $\mathrm{H}_{\bullet}(\Lambda M,M)$.

From now on let $\KK = \CC$ or $\KK = \HH$.
In case $\KK = \CC$ we set
$$   N = 2n, \quad \text{and}\quad  \lambda = 1    $$
and in case $\KK = \HH$ we set
$$   N = 4n, \quad \text{and}\quad  \lambda = 3 .      $$
We consider the free loop space of $M = \KK P^n$ where we consider $M$ as a symmetric space.
The symmetric Riemannian metric on $M$ induces a length functional $\mathcal{L}\colon \Lambda M\to \mathbb{R}$, see equation \eqref{eq_length_functional} and it is well-known that this is a Morse-Bott function on $\Lambda M$.
Moreover, the index and the nullity of all critical manifolds are finite, see \cite{ziller:1977}.
There is a group $G$ with a closed subgroup $K$ such that $M = G/K$ and $(G,K)$ is a Riemannian symmetric pair.
In particular the action of $G$ is a transitive action by isometries and $K$ is the isotropy group of a fixed basepoint $p_0\in M$.
The group $K$ then acts on $M$ by isometries as well and fixes the basepoint.
It is well-known that all geodesics on $M$ are closed and of the same prime length $l$.
Consequently, the critical values of the length functional are positive multiples of $l$.
If $k\in \NN$ then the critical set at level $a = kl$ is diffeomorphic to the unit tangent bundle $SM$, i.e.
$$  \Sigma^a \xrightarrow[]{\cong} S M, \qquad \gamma\mapsto \frac{\Dot{\gamma}(0)}{|\Dot{\gamma}(0)|}       $$
is a diffemorphism.
Moreover, Ziller argues in \cite{ziller:1977} that the group $G$ acts transitively on $\Sigma^a$ and this action equals the canonical action of $G$ on the unit tangent bundle.
Let $\gamma\in \Sigma^a$ be a closed geodesic at level $a = kl$ with $\gamma(0) = p_0$.
Then there is an underlying prime geodesic $\sigma \in \Sigma^l$ such that $\gamma = \sigma^k$.
Denote the isotropy group of $\gamma$ with respect to the action of $G$ on $\Sigma^a$ by $K_{\gamma}$.
In particular, this is a closed subgroup of $K$.
Then we have
$$    SM \cong \Sigma^a \cong G/K_{\gamma} .      $$
Furthermore, there is an induced action of the group $K$ on $\Sigma^a$ and the orbit of $\gamma$ is 
\begin{equation} \label{eq_k_oribt_sphere}
        K.\gamma \cong K/K_{\gamma} \cong \mathbb{S}^{N-1}    , 
\end{equation} 
which is the fiber of $SM$ over $p_0$.
In the following we will always write $\Sigma^k$ for the critical submanifold at level $kl$ instead of $\Sigma^{kl}$.
The index of $\Sigma^k$ is $$\mathrm{ind}(\Sigma^k) = k\lambda + (k-1)(N-1)  $$
see \cite[p. 167]{goresky:2009} and the nullity is equal to the dimension of $\Sigma^k$, i.e.
$$  \mathrm{null}(\Sigma^k) = 2N-1   . $$

On the figure-eight space $\Lambda M\times_M\Lambda M$ we consider the length function
$$  \mathcal{L}_2\colon \Lambda M\times_M \Lambda M\to [0,\infty),\quad \mathcal{L}_2(\eta_1,\eta_2) = \mathcal{L}(\eta_1) + \mathcal{L}(\eta_2)   $$
for $(\eta_1,\eta_2)\in\Lambda M\times_M\Lambda M$.
The critical manifolds in $\Lambda M\times_M \Lambda M$ are the sets of the form
$$  \Sigma_m\times_M \Sigma_{k-m}  = \Sigma_m\times \Sigma_{k-m} \cap \Lambda M\times_M \Lambda M , \quad k,m\in\mathbb{N}_0, \,\, k\geq m $$
and the critical values are again multiples of $l$.
Here the fiber product is taken with respect to the evaluation map at time $t=0$, i.e.
$$   \Sigma_m\times_M \Sigma_{k-m} = \{ (\gamma_1,\gamma_2)\in  \Sigma_m\times \Sigma_{k-m}\,|\, \gamma_1(0) = \gamma_2(0)\} .  $$
At level $kl$ the critical set is
$$    M\times_M \Sigma_k  \,\,  \sqcup \,\, \Sigma_1 \times_M \Sigma_{k-1} \,\, \sqcup \,\, \ldots\,\, \sqcup \,\, \Sigma_{k-1}\times_M  \Sigma_1\,\, \sqcup\,\, \Sigma_k\times_M M . $$
Note that by the relative construction  of the coproduct the components $M\times_M \Sigma_k$ and $\Sigma_k\times_M M$ will not show up in the course of the proof so we will not deal with them.
See also Remark \ref{remark_why_figure_eight}.
\begin{lemma} \label{lemma_l2}
    The length function $\mathcal{L}_2\colon \Lambda M\times_M \Lambda M\to [0,\infty)$ satisfies the Palais-Smale condition and is a Morse-Bott function.
    Moreover, we have 
    $$   \mathrm{ind}_{\Lambda\times_M \Lambda}(\eta_1,\eta_2) = \mathrm{ind}_{\Lambda}(\eta_1) + \mathrm{ind}_{\Lambda}(\eta_2) $$ and $$\mathrm{null}_{\Lambda\times_M\Lambda}(\eta_1,\eta_2) = \mathrm{null}_{\Lambda}(\eta_1) + \mathrm{null}(\eta_2) - N        $$
    for a critical point $(\eta_1,\eta_2)\in \Lambda M\times_M \Lambda M$ of the function $\mathcal{L}_2$.
\end{lemma}
\begin{proof}
    The function $$\mathcal{L}'\colon \Lambda \times \Lambda \to [0,\infty), \quad (\gamma_1,\gamma_2)\mapsto \mathcal{L}(\gamma_1)+\mathcal{L}(\gamma_2) $$
    clearly satisfies the Palais-Smale condition.
    Since $\Lambda M \times_M \Lambda M$ is a closed submanifold of $\Lambda \times \Lambda$ and since $\mathcal{L}_2$ is the restriction of $\widetilde{\mathcal{L}}$ it therefore follows that $\mathcal{L}_2$ also satisfies the Palais-Smale condition.

    In order to show that $\mathcal{L}_2$ is a Morse-Bott function, we need to show the following property.
    Let $(\eta_1,\eta_2)\in\Lambda M\times_M \Lambda M$ be a critical point of $\mathcal{L}_2$ and assume that it belong to the critical submanifold of the form $\Sigma^a\times_M \Sigma^b$ where $\Sigma^a$ and $\Sigma^b$ are critical submanifolds in $\Lambda M$ with respect to the Morse-Bott function $\mathcal{L}$.
    We need to show that the null space $T_{(\eta_1,\eta_2)}^0 \Lambda M\times_M \Lambda M$ is equal to the tangent space $T_{(\eta_1,\eta_2)}\Sigma^a \times_M \Sigma^b\subseteq T_{\eta_1,\eta_2} \Lambda M\times_M \Lambda M$.
    It is clear that we have
    $$      T_{(\eta_1,\eta_2)}\Sigma^a \times_M \Sigma^b\subseteq T_{\eta_1,\eta_2}^0 \Lambda M\times_M \Lambda M  . $$
    Arguing as in \cite[Section 2.5]{klingenberg:1995} one can see that $T_{(\eta_1,\eta_2)}^0\Lambda M\times_M \Lambda M$ can be characterized as 
    $$    T_{(\eta_1,\eta_2)}^0\Lambda \times_M \Lambda =  \{    (\xi_1,\xi_2)\in T_{\eta_1} \Lambda \oplus T_{\eta_2}\Lambda \,|\, \xi_1 (0) = \xi_2(0), \,\,\xi_1,\,\xi_2\,\,\text{periodic Jacobi fields}     \}  . $$
    Ziller shows in \cite[Section 2]{ziller:1977} that all periodic Jacobi fields along closed geodesics in a compact symmetric space are restrictions of Killing vector fields.
    Let $\xi_1,\xi_2$ be periodic Jacobi fields along $\eta_1$ and $\eta_2$, respectively.
    We assume without loss of generality that $\eta_1(0) =\eta_2(0) = p_0$ is the basepoint.
    If $\xi_1(0) = \xi_2(0) = 0$ then both Jacobi fields are restrictions of Killing fields on $M$ which are induced by the action of the group $K$.
    Since the action of $K\times K$ on $M\times M$ clearly preserves the diagonal $\Delta M$ it is clear that $(\xi_1,\xi_2)\in T_{(\eta_1,\eta_2)}\Sigma^a\times_M \Sigma^b$.
    If we have $\xi_1(0) = \xi_2(0) \neq 0$ then let us assume that $\nabla \xi_1(0) = \nabla \xi_2(0) = 0$.
    In this case one sees as in \cite[page 8]{ziller:1977} that both Jacobi fields are restrictions of the same Killing field, since the Killing field are determined by the element $\xi_1(0)$ in this case.
    Hence these Jacobi fields can be understood as restrictions of a Killing field of the diagonal group action $G\times M\times M\to M\times M$.
    Therefore in this case we also see that $(\xi_1,\xi_2)\in T_{(\eta_1,\eta_2)}\Sigma^a\times_M \Sigma^b$.
    Since the Jacobi fields of the above two types form a basis of $ T_{(\eta_1,\eta_2)}^0\Lambda \times_M \Lambda $ this shows the inclusion
    $$   T_{(\eta_1,\eta_2)}^0\Lambda \times_M \Lambda  \subseteq T_{(\eta_1,\eta_2)}\Sigma^a \times_M \Sigma^b .   $$
    Consequently, $\mathcal{L}_2$ is a Morse-Bott function and the claim for the nullity then follows directly from the dimensions of the critical submanifolds.
    Finally, for the indices, note that the index of a closed geodesic in a compact symmetric space is the same whether we consider it as a critical point in the based loop space or in the free loop space.
    Therefore, we get
    \begin{eqnarray*}
        \mathrm{ind}(\eta_1) + \mathrm{ind}(\eta_2) & = & \mathrm{ind}_{\Omega\times \Omega}((\eta_1,\eta_2)) \\ 
        & \leq & \mathrm{ind}_{\Lambda\times_M \Lambda}((\eta_1,\eta_2))  \\
        &\leq & \mathrm{ind}_{\Lambda\times\Lambda }((\eta_1,\eta_2))  \\
        &\leq & \mathrm{ind}(\eta_1) + \mathrm{ind}(\eta_2) 
    \end{eqnarray*}
    and thus we see that the inequalities are all equalities.
    This completes the proof.
\end{proof}

For $i,j\geq 1$ such that $i+j = k$ we have
$$   \Sigma_i\times_M \Sigma_j  \cong SM\times_M SM  = \{ (u,v)\in SM\times SM\,|\, \mathrm{pr}(u) = \mathrm{pr(v)} \}     $$
where $\mathrm{pr}\colon SM\to M$ is the canonical projection of the unit sphere bundle of $M$.
Moreover, the projection onto the first factor makes $SM\times_M SM$ into a sphere bundle over $SM$ which admits a global section $$SM\to SM\times_M SM, \quad u\mapsto (u,u),\quad u\in SM .$$
Since later on we shall use the cohomology ring of $SM\times_M SM$, we prove the following Lemma.

\begin{lemma} \label{lemma_cohomology_fiber_sm}
The rational cohomology ring of $SM\times_M SM$ is isomorphic to
$$  \mathrm{H}^{\bullet}(SM\times_M SM) \cong \frac{\QQ[{\alpha},{\beta},{\xi}]}{({\alpha}^n,{\beta}^2,{\xi}^2)}      $$
where $\mathrm{deg}({\alpha}) = \lambda + 1$, $\mathrm{deg}({\beta}) = N+\lambda$ and $\mathrm{deg}({\xi}) = N-1$.
\end{lemma}
\begin{proof}
From the Gysin sequence for $SM\to M$ we know that 
$$  \mathrm{H}^{\bullet}(SM) \cong \frac{\mathbb{Q}[\alpha,\beta]}{(\alpha^{n},\beta^2)}     $$
with $\mathrm{deg}(\alpha) = \lambda +1$ and $\mathrm{deg}(\beta) = N+\lambda$.
The manifold $SM\times_M SM$ is the total space of a sphere bundle over $SM$ with a global section, therefore it follows from the corresponding Gysin sequence that
$$ \mathrm{H}^{\bullet}(SM\times_M SM) \cong \mathrm{H}^{\bullet}(SM)\otimes \Lambda[{\xi}]    $$
with a generator ${\xi}$ of degree $\mathrm{deg}({\xi}) = N-1$.
This proves the claim.
\end{proof}
For $i\in\{0,\ldots,n-1\}$ we denote the homology class dual to the class ${\alpha}^i$ by $[{a}_i]$ and the dual of ${\alpha}^i{\beta}$ by $[{a}_i{b}]$.

Before we turn to the completing manifolds, let us note a property of the conjugate points along the closed geodesics in $M$.
With $\gamma$ and $\sigma$ as before, note that there is precisely one conjugate point $\sigma(\tfrac{1}{2}) = a$ along $\sigma$, see \cite[Proposition 3.35]{besse:1978}.
Moreover, the index of $\sigma$ is equal to $\lambda$, since the index of a closed geodesic on a compact symmetric space is equal to the sum of the multiplicity of the interior conjugate points, see again \cite[Proposition 3.35]{besse:1978} and \cite{ziller:1977}.
Denote the isotropy group of this point with respect to the action of $K$ by $K_a$.
It is well-known that
$$   \mathrm{dim}(K_a) > \mathrm{dim}(K_{\gamma})  ,    $$
see \cite[Theorem 4]{ziller:1977},
and that $\mathrm{dim}(K_a) - \mathrm{dim}(K_{\gamma})$ is equal to the index of $\sigma$ both as a geodesic loop in $\Omega M$ as well as a closed geodesic in $\Lambda M$.
\begin{lemma} \label{lemma_sphere_lambda}
    The homogeneous space $K_a/K_{\gamma}$ is diffeomorphic to the sphere $\mathbb{S}^{\lambda}$.
\end{lemma}
\begin{proof}
    By \cite[Proposition 3.35]{besse:1978} the set of first conjugate points along geodesics of the basepoint $p_0$ is equal to the cut locus of $p_0$.
    Moreover, we have
    $$  \mathrm{Cut}(\mathbb{K}P^n) \cong \mathbb{K}P^{n-1}  ,  $$
    see again \cite[Proposition 3.35]{besse:1978}.
    The set of tangent cut points 
    $$  S = \{v\in T_{p_0} M\,|\, \mathrm{exp}_{p_0}(v) \,\,\text{is the cut point of} \,\, t \mapsto \exp_p(tv) \}       $$
    is well-known to be the round sphere $\mathbb{S}^{N-1}$.
    Moreover, the exponential map induces a fiber bundle $\exp_{p_0}\colon \mathbb{S}^{N-1}\to \mathbb{K}P^{n-1}$, see \cite[Proposition 3.37]{besse:1978}.
    These are just the well-known fibrations of spheres over projective space, so it follows that the fiber is $\mathbb{S}^{\lambda}$.
    See also \cite[Theorem 5.29]{besse:1978} for details.
    We can understand these objects as homogeneous spaces, i.e. as we know from equation \eqref{eq_k_oribt_sphere} we have $K/K_{\gamma}\cong \mathbb{S}^{N-1}$ and it is clear that $K/K_a$ is diffeomorphic to the cut locus $\mathbb{K}P^{n-1}$.
    Therefore we see that that the fiber $K_{a}/K_{\gamma}$ of the fiber bundle $K/K_{\gamma}\to K/K_a$ is diffeomorphic to $\mathbb{S}^{\lambda}$.
\end{proof}

\section{Construction of Ziller's completing manifolds}\label{sec_construction_completing}

In this section we describe Ziller's completing manifolds, see \cite{ziller:1977}.
We describe the manifolds and the respective embeddings in detail.

Fix a closed geodesic $\gamma = \sigma^k\in \Lambda M$ of multiplicity $k$ starting at the basepoint $p_0\in M$.
Here, $\sigma$ is the underlying prime closed geodesic.
Consider the product
$$   W_k  = G\times K_a \times K\times  K_a \times K \ldots \times K_a       $$
with $2k$ factors in total. 
Throughout this section, we will follow the convention that the first element in the tuple 
$$  (g_0,x_1,\ldots,x_{2k-1})\in W_k      $$
is said to be at \textit{zero'th position}, the second one at \textit{first position} and so forth.
The element in zero'th position plays a special role since it lies in $G$, therefore we will denote it usually by $g_0$ while the other elements will be denoted by $x_i$.
There is a right action of the $2k$-fold product of $K_{\gamma}$ on $W_k$ given by
\begin{eqnarray*}    \chi : W_k \times K_{\gamma}^{2k} &\to& W_k      \\
                ( (g_0,x_1,\ldots,x_{2k-1}), (h_0,\ldots,h_{2k-1})) &\mapsto& (g_0 h_0,h_0^{-1} x_1 h_1,\ldots, h_{2k-2}^{-1} x_{2k-1} h_{2k-1})  .
\end{eqnarray*}
This action is free and proper and we consider the quotient space
$$  \Gamma_k = W_k/(K_{\gamma}^{2k})  .      $$
There is an embedding $f_k \colon \Gamma_k \to \Lambda M$ given by
$$   f_k([g_0,x_1,\ldots,x_{2k-1}])(t)  = \begin{cases} g_0. \gamma(t), & 0\leq t\leq \tfrac{1}{2k} 
\\
g_0x_1.\gamma(t), & \tfrac{1}{2k}\leq t \leq \tfrac{2}{2k}
\\
\,\,\,\,\, \vdots & \,\,\,\,\, \vdots
\\
g_0x_1\ldots x_{2k-1}.\gamma(t), & \tfrac{2k-1}{2k}\leq t \leq 1
\end{cases}      . $$
Note that the critical submanifold $\Sigma_k \cong G/K_{\gamma}$ can be seen as a submanifold of $\Gamma_k$ via the embedding
$$   s_{L,k}\colon G/K_{\gamma} \to \Gamma_k      ,\qquad   s_{L,k}([g]) = [ g,e,\ldots,e       ]    \,\,\,\text{for}\,\,\,[g]\in G/K_{\gamma}  .  $$
We define $L_k = s_{L,k}(G/K_{\gamma})$ and will identify $L_k$ and $G/K_{\gamma}$ in the following possibly without making the identification explicit.
There is a submersion 
\begin{equation} \label{eq_submersion_completing_mfld}
      p_{L,k} \colon \Gamma_k \to G/K_{\gamma}, \,\,\,\,p_{L,k}([g_0,x_1,\ldots,x_{2k-1}]) = [g_0]   \quad \text{for}\,\,\,    [g_0,x_1,\ldots,x_{2k-1}] \in \Gamma_k     
\end{equation}
and it is clear that $ p_{L,k}\circ s_{L,k} = \id_{G/K_{\gamma}}  $.
Moreover, we see that the composition $f_k\circ s_{L,k}$ is given by
$$    f_k \circ s_{L,k} ([g])(t) = g.\gamma(t)   \qquad \text{for}\,\,\,g\in G, \,\, t\in[0,1].    $$
Hence, the map $f_k \circ s_{L,k}$ is precisely the diffeomorphism $G/K_{\gamma}\cong \Sigma^k$.
Note that the only closed geodesics in the image of $f_k$ are precisely the closed geodesics in the critical submanifold $\Sigma_k$.
All other loops in $\mathrm{im}(f_k)$ are broken geodesics and hence they are not critical points of the energy functional.
Therefore, the flow of the energy functional decreases the value of $\mathcal{L}$ for all $\gamma\in \mathrm{im}(f_k)$ which are not in $\Sigma_k$.
Consequently, if we compose the embedding $f_k$ with an arbitrary short gradient flow of the length functional, we obtain a map of pairs
$$    (\Gamma_k, \Gamma_k \setminus L_k) \to (\Lambda M^{\leq kl},\Lambda M^{<kl})  .       $$

Since $$\mathrm{dim}(\Gamma_k) =  N +  k(N-1) + k (\mathrm{dim}(K_a)-\mathrm{dim}(K_{\gamma})) = \mathrm{ind}(\gamma) + 2N-1 = \mathrm{ind}(\gamma) + \mathrm{dim}(\Sigma_k)$$
we have shown that $\Gamma_k$ is a completing manifold for $\Sigma_k$ if we prove that it is orientable.
We will see the orientability later.

Remarkably, $\Gamma_k$ can also serve as a completing manifold in the Hilbert manifold 
$  \Lambda M\times_M \Lambda M$.
With $k$ and $\gamma$ as above, fix $1\leq m\leq k-1$ and consider the critical point $$(\gamma_1,\gamma_2) = (\sigma^m,\sigma^{k-m})\in \Lambda M\times_M \Lambda M,$$
where $\sigma$ is the underlying prime geodesic of $\gamma$.
The component of the critical set at level $k$ in $\Lambda M\times_M \Lambda M$ that contains $(\gamma_1,\gamma_2)$ is $\Sigma_m \times_M \Sigma_{k-m}$.
Note that
$$  \mathrm{ind}((\gamma_1,\gamma_2)) = \mathrm{ind}(\gamma_1) + \mathrm{ind}(\gamma_2) = k \, \mathrm{ind}(\sigma) + (k-2)(N-1) = \mathrm{ind}(\gamma) - (N-1) , $$
see Lemma \ref{lemma_l2}.
We want to see how $SM\times_M SM$ can be embedded into $\Gamma_k$.
Consider the right-action of $K_{\gamma}\times K_{\gamma}$ on $G\times K$ $$\chi' \colon (G\times K) \times (K_{\gamma}\times K_{\gamma}) \to  G\times K   $$
given by
$$   \chi'((g_0,x_{2m}),(h_0,h_{2m})) = (g_0 h_0,h_0^{-1} x_{2m} h_{2m})          $$
for $g_0\in G,x_{2m}\in K$ and $h_0,h_{2m}\in K_{\gamma}$.
Like the action $\chi$ above, this is a free and proper right action and we consider the quotient space
$   \mathcal{V} = G\times K/ \chi'      $.
\begin{lemma}
The manifold $\mathcal{V}$ is diffeomorphic to $SM\times_M SM$.
\end{lemma}
\begin{proof}
Recall that there is a diffeomorphism $ G/K_{\gamma} \to SM$ induced by the transitive action of $G$ on $SM$.
In particular we see that
\begin{equation} \label{eq_fiber_product_homogeneous}
      SM\times_M SM \cong \{ ([g_1],[g_2])\in G/K_{\gamma}\times G/K_{\gamma}\,|\, g_1^{-1} g_2\in K\} .       
\end{equation}
Moreover, let $E\subseteq G\times G$ be the submanifold
$$   E = \{ (g_1,g_2)\in G\times G\,|\, g_1^{-1} g_2\in K\} .     $$
It is clear that there is a submersion $E\to SM\times_M SM$ given by $(g_1,g_2) \mapsto ([g_1],[g_2])$.
Now, define maps
$$   \widetilde{\varphi}\colon G\times K\to E \quad \text{and}\quad \widetilde{\psi}\colon E\to G\times K   $$
by setting
$$   \widetilde{\varphi}(g,k) = (g,gk) \quad \text{for}\,\,\, g\in G,k\in K   $$
and
$$  \widetilde{\psi}(g_1,g_2) = (g_1,g_1^{-1}g_2)  \quad \text{for} \,\,\, (g_1,g_2)\in E.   $$
These maps factor through the submersions $G\times K\to \mathcal{V}$ and $E\to SM\times_M SM$ and therefore induce smooth maps
$$    \varphi\colon \mathcal{V}\to SM\times_M SM \quad \text{and}\quad \psi\colon SM\times_M SM\to \mathcal{V} .   $$
It is a direct computation that they are inverses of each other.
\end{proof}

Observe that there is an embedding
$$   s_{V,m} \colon \mathcal{V} \hookrightarrow \Gamma_k     $$
given by
$$   s_{V,m}( [g_0,x_{2m}] ) =  [g_0,e,\ldots,e,x_{2m},e,\ldots,e] \in \Gamma_k       $$
where $x_{2m}$ appears at the $2m$'th position.
We denote the image of $\mathcal{V}$ under this embedding by $V$.
There is a submersion $p_{V,m}\colon \Gamma_k \to \mathcal{V}$ given by
$$     p_{V,m} ([g_0,x_1,\ldots,x_{2k-1}]) = [g_0,x_1x_2\ldots x_{2m}] \in \mathcal{V}  .        $$
It is clear that $ p_{V,m}\circ s_{V,m} = \id_{\mathcal{V}} $.
We define a map $F_{k,m}\colon \Gamma_k \to \Lambda M\times_M \Lambda M$ as follows.
Let
\begin{equation} \label{eq_cut_k_cycle_map}
       F_{k,m}([g_0,x_1,\ldots,x_{2k-1}]) = (\eta_1,\eta_2)     \quad \text{for}\quad [g_0,x_1,\ldots,x_{2k-1}]\in \Gamma_k     
\end{equation}
where
$$  \eta_1(t) = \begin{cases}  g_0.\gamma_1(t ), & 0\leq t\leq \tfrac{1}{2m} 
\\
\vdots
\\
g_0 x_1\ldots x_{2m-1}.\gamma_1(t ), & \tfrac{2m-1}{2m}\leq t\leq 1
\end{cases}
$$
and
$$  
\eta_2(t) = \begin{cases}
g_0 x_1\ldots x_{2m-1}x_{2m}.\gamma_2(t) , & 0\leq t\leq\tfrac{1}{2k-2m}
\\
\vdots
\\
g_0 x_1\ldots x_{2k-1}.\gamma_2(t),& \tfrac{2k-2m-1}{2k-2m}\leq t\leq 1. 
\end{cases}
$$
It can be checked directly that $F_{k,m}$ is a continuous embedding.

\begin{lemma} \label{lemma_stragne_complete}
The embedding $F_{k,m} \colon  \Gamma_k\to \Lambda M\times_M \Lambda M$ maps $V$ homeomorphically onto the critical set $\Sigma_m \times_M \Sigma_{k-m}$. Moreover, the set of critical points in $\Lambda M\times_M \Lambda M$ in the image of $F_{k,m}$ is precisely the set $\Sigma_m\times_M \Sigma_{k-m}$.
\end{lemma} 
\begin{proof}
We show that the diagram
$$  
\begin{tikzcd}
\mathcal{V} \arrow[]{r}{\varphi} \arrow[]{d}{s_{V,m}} 
&  [3em]
SM\times_M SM  \arrow[]{d}{i_m}
\\
\Gamma_k \arrow[]{r}{F_{k,m} }
&
\Lambda M\times_M \Lambda M
\end{tikzcd}
$$
commutes, where the map $i_m$ is the inclusion of $SM\times_M SM$ into $\Lambda M\times_M \Lambda M$ as the critical set $\Sigma_m\times_M \Sigma_{k-m}$.
If we identify $SM\times_M SM$ with the subspace of $G/K_{\gamma}\times G/K_{\gamma}$ as in equation \eqref{eq_fiber_product_homogeneous} then
$$  i_m([g_1],[g_2]) = (g_1.\gamma_1,g_2.\gamma_2)\qquad\text{for} \,\,\,\, [g_1],[g_2]\in G/K_{\gamma}\,\,\,\text{with}\,\, g_1g_2^{-1}\in K  .     $$
With this identification the commutativity of the above diagram can be checked using the respective definitions.
The second statement can be checked from the definition of the map $F_{k,m}$.
\end{proof}

If we compose $F_{k,m}$ with an arbitrarily short flow of the gradient flow of the length functional $\mathcal{L}_2$ on $\Lambda M\times_M \Lambda M$ we obtain a map of pairs
$$    (\Gamma_k,\Gamma_k\setminus V) \to ((\Lambda M\times_M \Lambda M)^{\leq kl},(\Lambda M\times_M \Lambda M)^{<kl}) .     $$
Furthermore, we have
$$  \mathrm{dim}(\Gamma_k) = \mathrm{ind}(\gamma_1,\gamma_2) + \mathrm{dim}(\Sigma_m\times_M \Sigma_{k-m})    $$
so if we show that $\Gamma_k$ is orientable we see that $\Gamma_k$ is a completing manifold for $ \Sigma_m\times_M \Sigma_{k-m}$.

Note that in $\Lambda M\times_M \Lambda M$ we can also use the fiber product $\Gamma_m\times_M \Gamma_{k-m}$ as a completing manifold where the fiber product is taken with respect to the evaluation map
$$  \mathrm{ev}_k \colon \Gamma_k \to M,\quad [g_0,x_1,\ldots,x_{2k-1}] \mapsto g_0 K\in G/K\cong M     $$
for $1\leq l \leq k$.
Then one takes the map $$(f_m,f_{k-m})\colon \Gamma_m\times_M \Gamma_{k-m}\to \Lambda M\times_M \Lambda M$$
as an embedding of $\Gamma_m\times_M \Gamma_{k-m}$ and can check that this is again a completing manifold.
To conclude this section we want to show that the completing manifolds $\Gamma_k$ and $\Gamma_m\times_M \Gamma_{k-m}$ are equivalent.
Define a map $\Phi_m\colon \Gamma_k\to \Gamma_m\times_M \Gamma_{k-m}$ by
$$   \Phi_m ([g_0,x_1,\ldots, x_{2k-1}] ) = ( [g_0,x_1,\ldots,x_{2m-1}],[g_0 x_1\ldots x_{2m-1} x_{2m} , x_{2m+1},\ldots , x_{2k-1}])   $$
for $[g_0,x_1,\ldots,x_{2k-1}]\in \Gamma_k$.
Note that this is a well-defined and smooth map since it descends from an equivariant map $W_k\to W_m\times W_{k-m}$. 
Similarly, we define $\Psi_m\colon \Gamma_m\times_M \Gamma_{k-m}\to \Gamma_k$ by
\begin{eqnarray*}   & &  \Psi([g_0,x_1,\ldots ,x_{2m-1}],[g_{2m},\ldots ,x_{2k-1}]) \\ & &  =   [g_0,x_1,\ldots,x_{2m-1}, (g_0x_1\ldots x_{2m-1})^{-1}g_{2m},x_{2m+1},\ldots, x_{2k-1}]       \end{eqnarray*}
for $[g_0,x_1,\ldots ,x_{2m-1}]\in \Gamma_m, [g_{2m},\ldots,x_{2k-1}]\in\Gamma_{k-m}$ with $g_0^{-1}g_{2m}\in K$. 
One checks again that $\Psi_m$ is well-defined and smooth.

\begin{lemma} \label{lemma_equivalent_completing_manifolds}
The completing manifolds $\Gamma_k$ and $\Gamma_m\times_M \Gamma_{k-m}$ for the critical submanifold $\Sigma_m\times_M \Sigma_{k-m}$ are equivalent in the sense that the diagrams
$$
\begin{tikzcd}
SM\times_M SM \arrow[]{r}{} \arrow[dd , "\psi", "\cong"'] & \Gamma_m\times_M \Gamma_{k-m} \arrow[dd, "\Psi_m", "\cong"'] \arrow[]{dr}{(f_m,f_{k-m})}  & [2em]
\\
& & \Lambda M\times_M \Lambda M
\\
\mathcal{V} \arrow[]{r}{s_{V,m}} & \Gamma_k \arrow[swap]{ur}{F_{k,m}}  &
\end{tikzcd}
$$
and 
$$
\begin{tikzcd}
\Gamma_k \arrow[]{r}{\Phi_m} \arrow[swap]{d}{p_{V,m}} &  [2.5em] \Gamma_m\times_M \Gamma_{k-m} \arrow[]{d}{(p_{L,m},p_{L,k-m})} 
\\
\mathcal{V} \arrow[]{r}{\varphi} & SM\times_M SM
\end{tikzcd}
$$
commute.
Here, $p_{L,m}\colon \Gamma_m\to SM$ and $p_{L,k-m}\colon \Gamma_{k-m}\to SM$ are the submersions which are used for the completing manifold structure in the free loop space, see equation \eqref{eq_submersion_completing_mfld}.
In particular $\Phi$ is a diffeomorphism.
\end{lemma}
\begin{proof}
This can be checked directly by unwinding the definitions.
\end{proof}

\section{Cohomology of the completing manifolds}\label{sec_cohomology_completing}

The manifold $\Gamma_k$ is closely related to the $K$-cycles in the sense of Bott and Samelson, see \cite{bott:1958a}.
In this section we shall describe its homology and cohomology following the discussions of the cohomology of the $K$-cycles by Bott and Samelson \cite{bott:1958a} and by Araki \cite{araki:1962}.
Recall that $\Gamma_k$ is defined as the quotient of $W = (G\times K_a)\times (K\times K_a)^{k-1}$ modulo the action of $K_{\gamma}^{2k}$ via 
$$    \chi ((g_0,x_1,\ldots,x_{2k-1}), (h_0,\ldots,h_{2k-1})) = (g_0h_0,h_0^{-1}x_1 h_1^{-1},\ldots, h_{2k-2}^{-1}x_{2k-1} h_{2k-1})  .       $$
If we have any product of subgroups $K_i\subseteq G$, $i\in \{1,\ldots,m\}$ such that $K_{\gamma}\subseteq K_i$ for all $i\in \{1,\ldots,m\}$ there is an action of $(K_{\gamma})^m$ on this product given by
$$    ((x_1,\ldots,x_m), (h_1,\ldots,h_m)) \mapsto (x_1 h_1, h_1^{-1}x_2 h_2 ,\ldots , h_{m-1}^{-1} x_m h_m)      $$
for $(x_1,\ldots,x_m)\in K_1\times \ldots K_m$ and $(h_1,\ldots ,h_m)\in (K_{\gamma})^m$.
We denote the quotient by
$$   K_1 \times_{K_{\gamma}} K_2 \times_{K_{\gamma}} \ldots \times_{K_{\gamma}} (K_m/K_{\gamma}) .       $$
Note that we might have that the first group $K_1$ is the group $G$. All other groups will be subgroups of $K$.

\begin{lemma} \label{lemma_embed_sub_k_cycle}
Let $0\leq i_1 < i_2 < \ldots < i_m \leq 2k-1$ be integers and set $K_{0} = G$, $K_j = K$ for $j$ even, $j\geq 2$ and $K_j = K_a$ for $j$ odd.
Then the manifold 
$$   \Gamma^{i_1 i_2\ldots i_m} = K_{i_1}\times_{K_{\gamma}} K_{i_2}\times_{K_{\gamma}}\ldots\times_{K_{\gamma}} (K_{i_m}/K_{\gamma})       $$
can be embedded into $\Gamma_k$ via a map $s_{i_1\ldots i_m}\colon \Gamma^{i_1\ldots i_m}\hookrightarrow \Gamma_k$.
\end{lemma}
\begin{proof}
We define a map $\sigma_{i_1\ldots i_m}\colon K_{i_1}\times K_{i_2}\times \ldots K_{i_m}\to W_k$ by
$$   \sigma_{i_1\ldots i_m}(x_{i_1},\ldots x_{i_m})  = (e,\ldots,e,x_{i_1},e,\ldots, e,x_{i_2},e,\ldots, e,x_{i_m},e\ldots,e)          $$
where $x_{i_j}$ is at position $i_j$ for each $j\in\{1,\ldots,m\}$.
If $$(x_{i_1},\ldots ,x_{i_m})\in K_{i_1}\times \ldots \times K_{i_m} \quad \text{and}\quad (k_1,\ldots,k_m)\in K_{\gamma}^m$$ we have
\begin{eqnarray*}
  & & \sigma_{i_1\ldots i_m} ((x_{i_1} ,\ldots,x_{i_m}). (k_1,\ldots,k_m))  \\ &=& (e,\ldots, e,  x_{i_1} k_1, e,\ldots, e, k_1^{-1} x_{i_2} k_2,\ldots ) \\ &=&
 (e,\ldots, e, x_{i_1} k_1, k_1^{-1} k_1 ,\ldots, k_1^{-1} k_1, k_1^{-1} x_{i_2} k_2 ,k_2^{-1}k_2,\ldots )
 \\
 &=& 
\chi ( \sigma_{i_1\ldots i_m}(x_{i_1},\ldots,x_{i_m})  , (e,\ldots,e,k_1,k_1,\ldots,k_1,k_2,\ldots,k_2,\ldots)) .
 \end{eqnarray*}
Hence, $\sigma_{i_1\ldots i_m}$ is equivariant with respect to the action of $K_{\gamma}^m$ on $K_{i_1}\times\ldots \times K_{i_m}$ and the action $\chi$ of $K_{\gamma}^{2k}$ on $W_k$.
Therefore this yields a smooth map
\begin{equation} \label{eq_embedding_spheres}
      s_{i_1\ldots i_m}\colon \Gamma^{i_1\ldots i_m}\to \Gamma_k .   
\end{equation}
It is easy to check that this is an embedding.
\end{proof}

Let $P = (i_1,i_2,\ldots,i_m)$ with non-negative integers $0\leq i_1 < i_2 < \ldots < i_m \leq 2k-1$.
Then we say that $\Gamma^P = \Gamma^{i_1\ldots i_m}$ is a \textit{sub-$K$-cycle} of $\Gamma_k$.
From now on we will always identify a sub-$K$-cycle $\Gamma^P$ with its image in $\Gamma_k$ under the embeddings constructed in Lemma \ref{lemma_embed_sub_k_cycle}.
The manifold $\mathcal{V} \cong V$ which we defined above is an example of a sub-$K$-cycle.
If $P = (0,1,\ldots,m)$ for some $m\leq 2k-2$, then there are submersions
$$   p_{\Gamma,m} \colon \Gamma_k\to \Gamma^P      $$
given by
$$    p_{\Gamma,m}([g_0,x_1,\ldots,x_{2k-1}]) = [g_0,x_1,\ldots,x_m]         \quad \text{for}\,\,\, [g_0,x_1,\ldots,x_{2k-1}]\in \Gamma_k . $$
Hence, we get a chain of submersions
$$  \Gamma_k \to \Gamma^{0\ldots2k-2} \to \ldots \to \Gamma^{01} \to \Gamma^0 = G/K_{\gamma} .        $$
As shown in \cite[Theorem 2.4]{araki:1962} all these submersions are fiber bundles.
Moreover, each fiber bundle has a section which is given by the map
$$    \Gamma^{0\ldots m-1}\hookrightarrow \Gamma^{0\ldots m} , \quad  [g_0,x_1,\ldots,x_{m-1}] \mapsto [g_0,x_1,\ldots,x_{m-1},e] .      $$
\begin{definition}
Let 
$$E_0 \xrightarrow[]{\pi_0} E_1\xrightarrow[]{\pi_1} \ldots \xrightarrow[]{\pi_{m-1}} E_m = B $$
be a sequence of manifolds with each $\pi_i\colon E_i\to E_{i+1}$ being a sphere bundle.
Then we say that $E_0$ is an \textit{iterated sphere bundle over} $B$.
\end{definition}

\begin{prop}
The $K$-cycle $\Gamma_k$ is an iterated sphere bundle via the maps
$$  \Gamma_k \to \Gamma^{0\ldots2k-2} \to \ldots \to \Gamma^{01} \to \Gamma^0 = G/K_{\gamma} .        $$
\end{prop}
\begin{proof}
The fibers of the iterated fiber bundle are either $K/K_{\gamma}$ of $K_a/K_{\gamma}$.
We know that $K/K_{\gamma} \cong \mathbb{S}^{N-1}$ and $K_a/K_{\gamma} \cong \mathbb{S}^{\lambda}$, see Lemma \ref{lemma_sphere_lambda}.
Hence, $\Gamma_k$ is an iterated sphere bundle over $SM\cong G/K_{\gamma}$.
\end{proof}

In the situation of an iterated sphere bundle $E_0\to E_1\to\ldots\to E_m = B$, one can compute the cohomology of the total space $E_0$ by considering the Gysin sequences at each step.
In the following we consider homology and cohomology with rational coefficients.
We determine the cohomology ring of $\Gamma_k$.
Note that throughout the article we have fixed an orientation on $M$.
In particular this induces an orientation on the $\epsilon$-sphere around the basepoint $p_0$.
We shall denote the generator of its fundamental class by $[\mathbb{S}^{N-1}_{\epsilon}]$.
Furthermore, if $t\in I$, define $\mathrm{ev}_t\colon \Lambda M \to M$ to be the map $\mathrm{ev}_t(\gamma) = \gamma(t)$ for $\gamma\in \Lambda M$.

\begin{prop} \label{prop_cup_ring_gamma_k}
The cohomology ring of $\Gamma_k$ is isomorphic to
$$   \mathrm{H}^{\bullet}(\Gamma_k) \cong \frac{\mathbb{Q}[\alpha,\beta,\xi_1,\ldots,\xi_{2k-1}]}{ (\alpha^{n},\beta^2,\xi_1^2,\ldots,\xi_{2k-1}^2) }       $$
where $\mathrm{deg}(\alpha) = \lambda +1$, $\mathrm{deg}(\beta) = N+\lambda$, $\mathrm{deg}(\xi_{2i+1}) = \lambda$ for $i=0,\ldots,k-1$ and $\mathrm{deg}(\xi_{2i}) = N-1$ for $i=1,\ldots,k-1$.
In particular $\Gamma_k$ is orientable.

Furthermore, the class $\xi_{2i}$ and the dual class $[x_{2i}]$ in homology for $i\in\{1,\ldots,k-1\}$ can be chosen such that the following holds.
For $i\in\{1,\ldots,k-1\}$ one can choose a fundamental class $[\mathbb{S}^{N-1}]$ of $\mathbb{S}^{N-1}$ such that $[x_{2i}] = (s_{2i})_*[\mathbb{S}^{N-1}]$ and such that 
$$ (\mathrm{ev}_{t_i}|_{f_k(\Gamma_k)} \circ f_k \circ s_{2i})_*[\mathbb{S}^{N-1}] = [\mathbb{S}^{N-1}_{\epsilon}]  $$
where $\mathrm{ev}_{t_i}|_{f_k(\Gamma_k)}$ is understood as a map $\mathrm{ev}_{t_i}|_{f_k(\Gamma_k)}\colon f_k(\Gamma_k) \to B_{p_0}\setminus\{p_0\}$ for $t_i\in(\tfrac{i}{k} , \tfrac{i}{k} + \delta) $ for some small $\delta > 0$ and where $s_{2i}\colon \mathbb{S}^{N-1}\to \Gamma_k$ is the embedding constructed in Lemma \ref{lemma_embed_sub_k_cycle}.
\end{prop}
\begin{proof}
As we have seen before it follows from the Gysin sequence of $SM\to M$ that 
$$  \mathrm{H}^{\bullet}(SM) \cong \frac{\mathbb{Q}[\alpha,\beta]}{(\alpha^{n},\beta^2)}     $$
with $\mathrm{deg}(\alpha) = \lambda +1$ and $\mathrm{deg}(\beta) = N+\lambda$.
One can now determine the cohomology ring of $\Gamma_k$ by induction along the steps of the iterated sphere bundle.
Note that for a $k$-sphere bundle $E\to B$ with $k$ odd and which admits a global section one has
$$   \mathrm{H}^{\bullet}(E) \cong \mathrm{H}^{\bullet}(B) \otimes \Lambda_{\mathbb{Q}}[\xi]     $$
where $\mathrm{deg}(\xi) = k$.
Thus one iteratively obtains the cohomology ring.
For the orientations we note that at each step in the iterated sphere bundle we are free to choose the orientation of the new generator $\xi_i$.
As Araki argues using Gysin sequences, see \cite[Section 2]{araki:1962}, the homology class dual to the class $\xi_l$ for $l\in\{1,\ldots,2k-1\}$ can be described as follows.
Let $s_l\colon K_l/K_{\gamma}\hookrightarrow \Gamma_k$ be the embedding as in Lemma \ref{lemma_embed_sub_k_cycle}.
Then if we consider an orientation class $[K_l/K_{\gamma}]\in \mathrm{H}_{\bullet}(K/K_{\gamma})$ we have
$$  (s_l)_*[K_l/K_{\gamma}] = \pm [x_{l}]    .  $$
If $l=2i$ we have $K_{2i}/K_{\gamma}\cong \mathbb{S}^{N-1}$. 
Let $t_i\in (\tfrac{i}{k}\tfrac{i}{k}+\delta)$ with $\delta>0$ small then it can be seen directly that the map
$$  \mathrm{ev}_{t_i} \circ f_k\circ s_{2i}\colon   \mathbb{S}^{N-1} \to M    $$
maps $\mathbb{S}^{N-1}$ homeomorphically onto $\mathbb{S}^{N-1}_{\epsilon'}$ for some small $\epsilon'>0$.
Note that this property holds precisely because all geodesics in $M$ are closed and of the same prime length.

We can now choose an orientation class $[\mathbb{S}^{N-1}]$ and correspondingly the class $\xi_{2i}$ and its dual $[x_{2i}]$ such that
$$  [x_{2i}] = (s_{2i})_*[\mathbb{S}^{N-1}]\quad \text{and}\quad    (\mathrm{ev}_{t_i}\circ f_k\circ s_{2i})_*[\mathbb{S}^{N-1}] = [\mathbb{S}^{N-1}_{\epsilon'}]  $$
where we consider $\mathrm{ev}_{t_i}\circ f_k\circ s_{2i}$ as a map $\mathbb{S}^{N-1} \to \mathbb{S}^{N-1}_{\epsilon'}$.
\end{proof}

In the previous proposition we saw that the manifolds $\Gamma_k$, $k\in\mathbb{N}$ are orientable. 
This completes the proof that the $\Gamma_k$ are in fact completing manifolds.
We sum this up in the next corollary.

\begin{cor} Let $M = \KK P^n$ be a complex or quaternionic projective space and let $k,m\in\mathbb{N}$ with $k>  m$.
\begin{enumerate}
    \item 
The manifold $\Gamma_k$ with the embedding $f_k$ is a completing manifold for the critical set $\Sigma_k$ in $\Lambda M$.
\item
The manifold $\Gamma_k$ with the embedding $F_{k,m}$ is a completing manifold for the critical set $\Sigma_m\times_M \Sigma_{k-m}$ in $\Lambda M\times_M \Lambda M$. 
\end{enumerate}
\end{cor}
Note that the Corollary implies the perfectness of the Morse-Bott function $\mathcal{L}$ on $\Lambda M$.
This is one of the main results in \cite{ziller:1977}.

We now make the following notation convention for the generators in homology.
As shown above the classes
$$   \alpha^i \beta^j \xi_1^{l_1} \ldots \xi_{2k-1}^{l_{2k-1}} , \quad i\in\{0,\ldots,n-1\},\,\, j,l_i\in\{0,1\} \,\, \text{for}\,\,i\in\{1,\ldots,2k-1\}   $$
generate the cohomology of $\Gamma_k$ additively.
If $l_{i_1},\ldots, l_{i_p} = 1$ and $l_i = 0$ otherwise then we denote the dual of $ \alpha^i \beta^j \xi_1^{l_1} \ldots \xi_{2k-1}^{l_{2k-1}}$ in homology by
$$    [a_i  x_{i_1\ldots i_p}] \in \mathrm{H}_{\bullet}(\Gamma_k) \quad \text{if}\quad j=0 \quad \text{and} \quad  [a_i b x_{i_1\ldots i_p}]\quad \text{if}\quad j = 1.        $$
Recall that at level $kl$ the level homology $\mathrm{H}_{\bullet}(\Lambda M^{\leq kl},\Lambda M^{<kl})$ is isomorphic to the homology of the critical submanifold $L_k \cong SM$.
Moreover, as we have seen in Section \ref{sec_completing_manifolds} the map
$$   (p_{L,k})_! \colon \mathrm{H}_{i-\lambda_k}(L_k) \to \mathrm{H}_i(\Gamma_k) 
$$
is injective, where $$\lambda_k = \mathrm{ind}(\gamma)    = k\lambda + (k-1)(N-1)  .  $$ 
Hence, we obtain the generators of $\mathrm{H}_{\bullet}(\Lambda M)$ which come from level $k$ by considering the map $(p_{L,k})_!$.
Recall that $L_k \cong SM$ and as we have seen its cohomology ring is
$$   \mathrm{H}^{\bullet}(L_k;\QQ) \cong \frac{\QQ[{\alpha},{\beta}]}{({\alpha}^n,{\beta}^2)}    $$
with $\mathrm{deg}({\alpha}) = \lambda + 1$ and $\mathrm{deg}({\beta}) = N + \lambda$.
In homology we choose dual generators and denote them by $[{a}_0],\ldots, [{a}_{n-1}]\in \mathrm{H}_{\bullet}(SM;\QQ)$ with $\mathrm{deg}([{a_i}]) = (\lambda +1) i$ and $[{a}_0{b}],\ldots,[{a}_{n-1}{b}] \in \mathrm{H}_{\bullet}(SM;\QQ)$ with $\mathrm{deg}([{a_i}{b}]) = N +\lambda + i(\lambda + 1)$.
In particular, we can choose these generators such that under the embedding
$$  s_{L,k}\colon G/K_{\gamma} \cong SM \hookrightarrow \Gamma_k      $$
we have
$$   (s_{L,k})_* [{a}_i] = [ a_i] \quad \text{and}\quad (s_{L,k})_*[{a}_i{b}] = [a_i b]      $$
and in cohomology
\begin{equation} \label{eq_pullback_generators}
       (p_{L,k})^*{\alpha}^i = \alpha^i \quad \text{and} \quad (p_{L,k})^* {\beta} = \beta 
\end{equation}
where $p_{L,k}\colon \Gamma_k\to L_k$ is the retraction.
The above formulas also justify the misuse of notation, since e.g. the cohomology class $\alpha$ has a double meaning, but as $(p_{L,k})_*\colon \mathrm{H}^{\bullet}(L_k)\to \mathrm{H}^{\bullet}(\Gamma_k)$ is injective, it is reasonable to identify $\alpha\in\mathrm{H}^{\bullet}(L_k)$ with its image under this injection.
We choose the orientation of $L_k$ by choosing the class $[{a}_{n-1}{b}]\in\mathrm{H}_{2N-1}(L_k)$ as fundamental class for all $k\in\mathbb{N}$ and as fundamental class for $\Gamma_k$ we choose the class $$[a_{n-1}b x_{1\ldots 2k-1}]\in \mathrm{H}_{2N-1 + \lambda_k}(\Gamma_k)  . $$
Recall that  the Gysin map $(p_{L,k})_!\colon \mathrm{H}_{\bullet}(L_k)\to \mathrm{H}_{\bullet+\lambda}(\Gamma_k)$ is defined as the composition 
$$(p_{L,k})_! \colon \mathrm{H }_j(L_k) \xrightarrow[]{(PD_{L_k})^{-1}}  \mathrm{H}^{\mathrm{dim}(L_k) - j}(L_k) \xrightarrow[]{p_{L,k}^*} \mathrm{H}^{\mathrm{dim}(L_k)-j}(\Gamma_k) \xrightarrow[]{PD_{\Gamma_k}} \mathrm{H}_{j + \lambda}(\Gamma_k) . $$
Using Proposition \ref{prop_cup_ring_gamma_k} we can now compute the map $(p_{L,k})_!$.
\begin{lemma}
With the above notation the following equations hold
$$  (p_{L,k})_!([{a}_i]) = - [a_ix_{1\ldots 2k-1}] \in \mathrm{H}_{\lambda_k + i(\lambda +1)}(\Gamma_k)       $$
and 
$$     (p_{L,k})_!([{a}_i{b}]) =  [a_i b x_{1\ldots 2k-1}] \in \mathrm{H}_{\lambda_k + (i+1)(\lambda +1) + N-1 }(\Gamma_k)   . $$
\end{lemma}
\begin{proof}
    We just consider the first case, the second case is analogous.
    The Poincaré dual of $[{a}_i]\in \mathrm{H}_{\bullet}(L_k)$ is the cohomology class ${\alpha}^{n-1-i}{\beta}\in \mathrm{H}^{\bullet}(L_k)$.
    By equation \eqref{eq_pullback_generators} we have 
    $$   (p_{L,k})^* ({\alpha}^{n-1-i}{\beta})  = \alpha^{n-1-i}\beta \in \mathrm{H}^{\bullet}(\Gamma_k) .      $$
    We now need to compute the Poincaré dual $X = PD_{\Gamma_k}(\alpha^{n-1-i}\beta)$ of this latter class.
    We compute the Kronecker pairing
    \begin{eqnarray*}
        \langle \alpha^i \xi_1 \ldots \xi_{2k-1}, X \rangle & = & \langle \alpha^i\xi_1\ldots \xi_{2k-1} \,,\, \alpha^{n-1-i}\beta \cap [\Gamma_k] \rangle \\
        &=& 
        \langle \alpha^i \xi_1\ldots \xi_{2k-1} \, \cup \, \alpha^{n-1-i}\beta \,,\,  [\Gamma_k]\rangle .
    \end{eqnarray*}
    Now, since $\alpha$ is of even degree and $\beta$ is of odd degree we get by graded commutativity
    $$   \alpha^i \xi_1\ldots \xi_{2k-1} \, \cup \, \alpha^{n-1-i}\beta =  - \alpha^{n-1}\beta \xi_1\ldots \xi_{2k-1}   $$
    and therefore
     $$\langle \alpha^i \xi_1 \ldots \xi_{2k-1}, X \rangle = -1    .$$
     It follows that $X = - [a_i b x_{1\ldots 2k-1}]$.
\end{proof}

We define classes
$$   A_k^i = (f_k)_* \big(    - [a_ix_{1\ldots 2k-1}]   \big)   \in \mathrm{H}_{\lambda_k + i(\lambda +1)}(\Lambda M)    $$
and 
$$   B_k^i = (f_k)_* \big(    [a_i b x_{1\ldots 2k-1}]     \big) \in \mathrm{H}_{\lambda_k + (i+1)(\lambda +1) + N-1 }(\Lambda M)     $$
for $k\in \NN$ and $i\in\{0,\ldots,n-1\}$.
Note that the degree of all $A^i_k$ is odd, while the degree of all $B_k^i$ is even.
The following is then clear by the construction of the completing manifolds.
\begin{prop} \label{prop_generators_free_loop_space}
The homology of the free loop space relative to the constant loops is generated by the image of the set
$$ \{ A_k^i \in \mathrm{H}_{\bullet}(\Lambda M) \,|\,k\in \NN, i\in\{0,\ldots,n-1\} \} \,\,\cup\,\,  \{ B_k^i\in \mathrm{H}_{\bullet}(\Lambda M) \,|\,k\in \NN, i\in\{0,\ldots,n-1\} \}       $$
in the relative homology $\mathrm{H}_{\bullet}(\Lambda M,M)$.
\end{prop}

\section{Computation of the coproduct} \label{sec_computation}

In this section we explicitly compute the string topology coproduct.
As we have seen in the previous section we can explicitly describe a set of generators of the homology $\mathrm{H}_{\bullet}(\Lambda M, M)$ via the completing manifolds $\Gamma_k$, $k\in \NN$.
We therefore express all the steps in the definition of the coproduct in intrinsic terms of the manifolds $\Gamma_k$.

First, we want to pull back the class $\tau_{\Lambda}\in \mathrm{H}^N(U_{\Lambda},U_{\Lambda,\geq\epsilon_0})$ via $f_k$ to a class which can be described in terms of the cohomology of $\Gamma_k$.
We consider the preimage $(f_k,\id_I)^{-1}(U_{\Lambda})\subseteq \Gamma_k\times I$.
This set can be described explicitly as follows.
Since all loops in the image $f_k(\Gamma_k)$ are broken geodesics, there is a small $\delta>0$ such that
$$  (f_k,\id_I)^{-1}(U_{\Lambda}) = \Gamma_k \times \Big( [0,\delta) \cup (\tfrac{2}{2k}-\delta,\tfrac{2}{2k}+\delta) \cup \ldots \cup (\tfrac{2k-2}{2k}-\delta,\tfrac{2k-2}{2k}+\delta)\cup (1-\delta,1]\Big)  .          $$
This is because every other conjugate point on a closed geodesic starting at the basepoint is the basepoint itself.
Clearly, $\delta>0$ is so small that the open intervals are disjoint.
Similarly, there is a $\delta_0>0$ with $\delta_0<\delta$ such that
\begin{eqnarray*}
  (f_k,\id_I)^{-1}(U_{\Lambda,\geq\epsilon_0}) = \Gamma_k \times   \Big( [\delta_0,\delta) \cup (\tfrac{2}{2k}-\delta, \tfrac{2}{2k}-\delta_0]\cup [\tfrac{2}{2k}+\delta_0,\tfrac{2}{2k} +\delta) \cup   \\
  \ldots \cup (\tfrac{2k-2}{2k}-\delta,\tfrac{2k-2}{2k}-\delta_0]\cup [ \tfrac{2k-2}{2k}+\delta_0,\tfrac{2k-2}{2k}+\delta) \cup (1-\delta,1-\delta_0] \Big)  .
\end{eqnarray*}
To make the bookkeeping easier, let us define $$I_m = (\tfrac{2m}{2k}-\delta,\tfrac{2m}{2k}+\delta)$$ and $$J_m = (\tfrac{2m}{2k} - \delta,\tfrac{2m}{2k}-\delta_0] \cup [\tfrac{2m}{2k}+\delta_0,\tfrac{2m}{2k}+\delta)$$ for $m = 1,\ldots,k-1$.
We set $$U_{\Gamma_k} = (f_k,\id_I)^{-1}(U_{\Lambda})\quad  \text{ and} \quad U_{\Gamma_k,\geq\epsilon_0} =  (f_k,\id_I)^{-1}(U_{\Lambda,\geq\epsilon_0}), $$ then we have
$$  (U_{\Gamma_k},U_{\Gamma_k,\geq\epsilon_0}) = \Gamma_k \times  \Big( ([0,\delta),[\delta_0,\delta)) \sqcup \bigsqcup_{m=1}^{k-1} (I_m,J_m) \sqcup ((1-\delta,1],(1-\delta,1-\delta_0] ) \Big)  .         $$
The pairs 
$$   ([0,\delta),[\delta_0,\delta)) \quad \text{and}\quad ((1-\delta,1],(1-\delta,1-\delta_0])     $$
have trivial homology for obvious reasons.
For $m\in\{1,\ldots,k-1\}$ we have
$$  H_i(I_m,J_m) \cong \begin{cases}  \QQ, &    i=1\\ 0& \text{else} .    \end{cases}          $$
We choose positively oriented generators $[I_m]$ of $\mathrm{H}_1(I_m,J_m)$ and dual cohomology classes
$$  \eta_m\in \mathrm{H}^1(I_m, J_m) .     $$
We now want express the cohomology class $\tau_k =  (f_k,\id_I)^{*}\tau_{\Lambda} \in \mathrm{H}^{N}(U_{\Gamma_k},U_{\Gamma_k,\geq\epsilon_0})$ in terms of the cohomology of $\Gamma_k$.
Recall that the cap product which we use in the definition of the string topology coproduct is a particular relative version of the ordinary cap product.
We refer to Appendix \ref{appendix_cap} for details, see also \cite[Appendix A]{hingston:2017}.

\begin{lemma} \label{lemma_thom_class_k_cycle}
The pullback of the class $\tau_{\Lambda}\in \mathrm{H}^{N}(U_{\Lambda},U_{\Lambda,\geq\epsilon_0})$ under the map 
$$   (f_k,\id_I)\colon (U_{\Gamma_k},U_{\Gamma_k,\geq\epsilon_0}) \to (U_{\Lambda},U_{\Lambda,\geq\epsilon_0})     $$
is given by 
$$  \tau_k =  (f_k,\id_I)^*\tau_\Lambda = \sum_{m=1}^{k-1}   \xi_{2m}\times \eta_m .      $$
\end{lemma}

This key lemma is proved in Appendix \ref{appendix_thom_class}.
In the following let $Y$ be one of the classes
$$   A_k^i = (f_k)_* \big(    - [a_ix_{1\ldots 2k-1}]   \big)   \in \mathrm{H}_{\bullet}(\Lambda M)    $$
or
$$   B_k^i = (f_k)_* \big(    [a_i b x_{1\ldots 2k-1}]     \big) \in \mathrm{H}_{\bullet}(\Lambda M)     $$
for $k\in \NN$, $i\in\{0,\ldots,n-1\}$.
We write
$$   Y = (f_k)_* X     $$
where $X\in \mathrm{H}_{\bullet}(\Gamma_k)$ is the respective homology class in $\Gamma_k$.
In order to compute $\cpr Y$ we need to consider
\begin{eqnarray*}
  \tau_{\Lambda} \cap (Y\times [I])  &=& (f_k,\id_I)_*( (f_k,\id_I)^*\tau_{\Lambda} \cap (X\times [I])) \\ &=& (f_k,\id_I)_*(\tau_k\cap (X\times [I]) )    
\end{eqnarray*}
by naturality of the cap product, see Proposition \ref{prop_nat_cap}.
As seen in Lemma \ref{lemma_thom_class_k_cycle}, we obtain
$$    \tau_k \cap (X\times [I]) =    \sum_{m=1}^{k-1}  \big( (\xi_{2m}\times \eta_m) \cap (X\times [I])       \big) .  $$
By the compatibility of the cross and the cap product, see Proposition \ref{prop_compat_cross_cap}, we have
$$   (\xi_{2m}\times \eta_m)\cap (X\times [I]) = (-1)^{\mathrm{deg}(X)} (\xi_{2m}\cap X)\times (\eta_m\cap [I]) .      $$
By the construction of the relative cap product we see  that
$$   \eta_m\cap [I] = [t_m]\in \mathrm{H}_0(I_m)      $$
where $[t_m]$ is a generator of $\mathrm{H}_0(I_m)$, see Example \ref{ex_relative_cap_interval}.
\begin{lemma} \label{lemma_cap_thom_explicit}
The relative cap product yields
$$   (\xi_{2m}\times \eta_m)\cap \big (  - [a_i x_{1\ldots 2k-1}]    \times [I] \big) = - [a_i x_{1\ldots 2m-1\, 2m+1\ldots 2k-1}]   \times [t_m]     $$
and 
$$   (\xi_{2m}\times \eta_m)\cap \big(     [a_i b x_{1\ldots 2k-1}] ) \times [I] \big) = - [a_i b x_{1\ldots 2m-1\, 2m+1\ldots 2k-1}]  \times [t_m]      $$
\end{lemma}
\begin{proof}
Using Proposition \ref{prop_cup_ring_gamma_k} we have
$$  \xi_{2m} \cap \big(   - [a_i x_{1\ldots 2k-1}]    \big) = [a_i x_{1\ldots 2m-1 \, 2m+1\ldots 2k-1}]    $$
and 
$$  \xi_{2m} \cap     [a_i b x_{1\ldots 2k-1}]    = - [a_i b x_{1\ldots 2m-1 \, 2m+1\ldots 2k-1}]    $$
and this yields the claim.
\end{proof}
For convenience of notation we shall write
$$   (\xi_{2m}\times \eta_m)\cap (X\times [I])  = X_m\times [t_m]   $$
and plug in the respective classes later using the above Lemma.
Then we have
$$   \tau_k\cap(X\times [I]) = \sum_{m=1}^{k-1} X_m\times [ t_m] .      $$
Fix an $m\in\{1,\ldots,k-1\}$.
To finish the computation of $\cpr Y$, we need to determine the effect of the retraction map $\mathrm{R}_{GH}$ and of the cutting map on $X_m\times [t_m] $.
First note that the diagram
$$   
\begin{tikzcd}
\mathrm{H}_{\bullet}(\Gamma_k\times I_m) \arrow[]{r}{(f_k,\id_I)_*} 
\arrow[swap]{d}{(\id_{\Gamma_k},\sigma_m)_*}
& 
\mathrm{H}_{\bullet}(U_{\Lambda}) \arrow[]{d}{=} \arrow[]{r}{}
&
\mathrm{H}_{\bullet}(U_{\Lambda},M\times I\cup \Lambda\times\partial I)
\arrow[]{d}{=}
\\
\mathrm{H}_{\bullet}(\Gamma_k\times \{\tfrac{2m}{2k}\}) \arrow[]{r}{(f_k,\id_I)_*}
\arrow[swap]{d}{ =}
&
\mathrm{H}_{\bullet}(U_{\Lambda}) \arrow[]{r}{} \arrow[]{d}{(R_{GH})_*}
&
\mathrm{H}_{\bullet}(U_{\Lambda},M\times I \cup \Lambda\times\partial I)
\arrow[]{d}{(\mathrm{R}_{GH})_*}
\\
\mathrm{H}_{\bullet}(\Gamma_k\times \{\tfrac{2m}{2k}\}) \arrow[]{r}{(f_k,\id_I)_*}
&
\mathrm{H}_{\bullet}(F_{\Lambda}) \arrow[]{r}{}
&
\mathrm{H}_{\bullet}(F_{\Lambda}, M \times I \cup \Lambda\times\partial I)
\end{tikzcd}
$$
commutes, where $\sigma_m\colon I_m\to \{\tfrac{2m}{2k}\}$ is the constant map.
To complete the computation, we need to characterize
$$   \mathrm{cut}_*(f_k,\id_I)_* (X_m\times [\tfrac{2m}{2k}]) \in \mathrm{H}_{\bullet}(\Lambda M \times \Lambda M,\Lambda M\times M\cup M\times\Lambda M) . $$

A direct computation shows that the map
$$  \widetilde{\mathrm{cut}} \circ f_k \, \colon \,\, \Gamma_k\times \{\tfrac{2m}{2k}\} \to \Lambda M\times_M \Lambda M      $$
is equal to the map $F_{k,m}\colon \Gamma_k \to \Lambda M\times_M \Lambda M$, which was defined in equation \eqref{eq_cut_k_cycle_map}, up to the obvious identification $\Gamma_k\cong \Gamma_k\times \{\tfrac{2m}{2k}\}$.
This shows that 
$$     \mathrm{\widetilde{cut}}_*(f_k,\id_I)_* ( X_m\times [\tfrac{2m}{2k}]) = (F_{k,m})_* X_m .         $$
Here and in the following $[\tfrac{2m}{2k}]\in\mathrm{H}_0(\{\tfrac{2m}{2k}\})$ denotes the canonical generator.
We now want to express the class $(\iota \circ F_{k,m})_* X_m \in \mathrm{H}_{\bullet}(\Lambda M\times \Lambda M)$ as a product of the generators $A_k^i$ and $B_k^i$.
Here, $\iota\colon \Lambda M\times_M \Lambda M\hookrightarrow \Lambda M\times \Lambda M$ is the inclusion of the figure-eight space.
In order to do so we need the following lemma.

\begin{lemma}\label{lemma_comm_diag_gysin}
The following diagram commutes
$$   
\begin{tikzcd}
\mathrm{H}_{i-\lambda_k +(N-1)}(V) \arrow[]{r}{(p_{V,m})_!} \arrow[d, "\varphi_*", "\cong"']
& [2.5em]
\mathrm{H}_i(\Gamma_k) \arrow[d, "\Phi_*", "\cong"'] \arrow[]{rd}{(F_{k,m})_*} & [2.5em]
\\
\mathrm{H}_{i-\lambda_k + (N-1)}(SM\times_M SM) \arrow[swap]{r}{(p_{L,m},p_{L,k-m})_!}  \arrow[]{d}{}
&
\mathrm{H}_i(\Gamma_m\times_M \Gamma_{k-m}) \arrow[swap]{r}{(f_m,f_{k-m})_*}  \arrow[]{d}{}  & 
\mathrm{H}_i(\Lambda M\times_M \Lambda M) \arrow[]{d}{}
\\
\mathrm{H}_{i-\lambda_k + (N-1)}(SM\times SM) \arrow[swap]{r}{(p_{L,m},p_{L,k-m})_!} &
\mathrm{H}_{i}(\Gamma_m\times \Gamma_{k-m}) \arrow[swap]{r}{(f_m,f_{k-m})_*} 
&
\mathrm{H}_i(\Lambda M\times \Lambda M)
\end{tikzcd}
$$
where the vertical arrows in the lower row are induced by the respective inclusions.
\end{lemma}
\begin{proof}
    The commutativity of all subdiagrams is clear apart from the lower left square.
    In order to show that the lower left square commutes, we first consider the diagram
    $$
    \begin{tikzcd}
        \mathrm{H}_{i-\lambda_k + (N-1)}(SM\times_M SM)   \arrow[]{d}{} 
        & [2.0em]
        \mathrm{H}_i(\Gamma_m\times_M \Gamma_{k-m})  \arrow[]{d}{} \arrow[swap]{l}{\mathrm{Th}'}
        \\
        \mathrm{H}_{i-\lambda_k + (N-1)}(SM\times SM) 
        &
        \mathrm{H}_i(\Gamma_m\times \Gamma_{k-m}) \arrow[swap]{l}{\mathrm{Th}}
    \end{tikzcd}
    $$
    where $\mathrm{Th}\colon \mathrm{H}_i(\Gamma_m\times \Gamma_{k-m})$ is the map
    \begin{eqnarray*}
            \mathrm{H}_i(\Gamma_m\times \Gamma_{k-m}) &\xrightarrow[]{\hphantom{excision}} & \mathrm{H}_i(\Gamma_m\times \Gamma_{k-m}, \Gamma_m\times \Gamma_{k-m}\setminus SM\times SM) \\ & \xrightarrow[]{\text{excision}} & \mathrm{H}_i(U,U\setminus SM\times SM) \\ & \xrightarrow[]{\hphantom{b}\text{Thom}\hphantom{b}} & \mathrm{H}_{i-\lambda_k + (N-1)}(SM\times SM) .      
    \end{eqnarray*}
    Here, $U$ is a tubular neighborhood of $SM\times SM$ in $\Gamma_m\times \Gamma_{k-m}$ and $$\mathrm{Thom}\colon \mathrm{H}_{\bullet}(U,U\setminus SM\times SM) \to \mathrm{H}_{\bullet-\lambda_k +(N-1)}(SM\times SM)$$ is the Thom isomorphism.
    The map $\mathrm{Th}'$ is defined analogously.
    In particular, we note that the normal bundle of $SM\times_M SM \hookrightarrow \Gamma_m\times_M \Gamma_{k-m}$ is the pullback of the normal bundle of $SM\times SM\hookrightarrow \Gamma_m\times \Gamma_{k-m}$ along the inclusion $SM \times_M SM \hookrightarrow SM\times SM$.
    Therefore the above diagram commutes.
    Now, note that the map $\mathrm{Th}$ agrees with the Gysin map $(s_m,s_{k-m})_!$.
    This follows from \cite[Theorem VI.11.3]{bredon:2013}.
    Note that in this reference it is only claimed that the two maps agree up to sign, but one can determine the sign from the proof.
    Applied to our present case the sign is $(-1)^{c (\lambda_k - (N-1))}$ for some integer $c\in \mathbb{Z}$.
    Recall that the index $\lambda_k$ is odd for all $k\in\mathbb{N}$.
    Consequently, the codimension 
    $$   \mathrm{codim}(SM\times SM\hookrightarrow \Gamma_m\times \Gamma_{k-m}) = \lambda_k - (N-1)    $$
    is even for all $k\in\mathbb{N}$, so we see that the sign $(-1)^{c (\lambda_k - (N-1))}$ is even and thus the maps $\mathrm{Th}$ and $(s_{L,m},s_{L,k-m})_!$ agree.
    Now, let $Z\in\mathrm{H}_{\bullet}(SM\times_M SM)$.
    Then we have
    $$     Z =   ((p_{L,m},p_{L,k-m}) \, \circ \, (s_{L,m},s_{L,k-m}) )_! \,Z =    (s_{L,m},s_{L,k-m})_! \circ (p_{L,m},p_{L,k-m})_! \, Z . $$
    Moreover, let us denote the inclusion $SM\times_M SM\hookrightarrow SM\times SM$ by $i_1$ and denote the inclusion $\Gamma_m\times_M \Gamma_{k-m}\hookrightarrow \Gamma_m\times \Gamma_{k-m}$ by $i_2$.
    Then we get
    \begin{eqnarray*}
            (p_{L,m},p_{L,k-m})_! \circ (i_1)_* Z &=&  (p_{L,m},p_{L,k-m})_! \circ (i_1)_*  (s_{L,m},s_{L,k-m})_! \circ (p_{L,m},p_{L,k-m})_! Z  \\ &=&     (p_{L,m},p_{L,k-m})_! \circ  (s_{L,m},s_{L,k-m})_! \circ (i_2)_* \circ (p_{L,m},p_{L,k-m})_! Z \\&=& (i_2)_* \circ (p_{L,m},p_{L,k-m})_! Z  
    \end{eqnarray*} 
    and this shows the commutativity of the lower left square.    
\end{proof}

Recall that $\Gamma_k$ together with the embedding $F_{k,m}$ is a completing manifold for the critical set $\Sigma_m\times_M \Sigma_{k-m}\cong V$.
In Lemma \ref{lemma_cohomology_fiber_sm} we determined the cohomology ring of $V$.
We make the following choice of orientation.
The generator $\xi\in \mathrm{H}^{N-1}(V)$ is chosen in such a way that it pulls back to the generator $\xi_{2m}$ under the map $p_{V,m}$.
Note that by the choice of orientations for the classes $\xi_{2i}$, $i\in\{1,\ldots,k-1\}$ in Proposition \ref{prop_cup_ring_gamma_k} this is a consistent choice.
We choose the orientation for $V$ such that $\alpha^{n-1}\beta\xi$ is a fundamental cohomology class.

\begin{lemma} \label{lemma_gysin_v}
Using the notation for the generators of $\mathrm{H}_{\bullet}(SM\times_M SM)$ as in the paragraph following Lemma \ref{lemma_cohomology_fiber_sm} we have
$$   (p_{V,m})_! [ {a}_i] = - [a_ix_{1\ldots 2m-1\, 2m+1 \ldots 2k-1}]        $$
and
$$   (p_{V,m})_![{a}_i{b}] =  - [a_i b x_{1\ldots 2m-1 \, 2m+1 \ldots 2k-1}] .     $$
\end{lemma}
\begin{proof}
We only consider the first case, the second one is analogous.
The Poincaré dual of $[{a}_i]$ is the class $\alpha^{n-1-i}\beta\xi\in\mathrm{H}_{\bullet}(V)$.
Moreover, by our orientation convention, the pullback of this cohomology class is
$$    (p_{V,m})^* (\alpha^{n-1-i}\beta\xi) = \alpha^{n-1-i}\beta\xi_{2m} .    $$
Now we need to compute the Poincaré dual of this class in $\Gamma_{k}$, i.e.
$$    PD_{\Gamma_k}( \alpha^{n-1-i}\beta\xi_{2m}) = \alpha^{n-1-i}\beta\xi_{2m} \cap [\Gamma_k] .   $$
By using the graded commutativity of the cup product, we see that 
$$     \langle \alpha^i \xi_1\ldots \xi_{2m-1}\xi_{2m+1}\ldots \xi_{2k-1}, \alpha^{n-1-1}\beta\xi_{2m} \cap [\Gamma_k]\rangle = \langle -\alpha^{n-1}\beta\xi_1\ldots \xi_{2k-1},[\Gamma_k]\rangle = -1 .      $$
Therefore we see that
$$   \alpha^{n-1-1}\beta\xi_{2m} \cap [\Gamma_k] =  - [a_i x_{1\ldots 2m-1 \, 2m+1\ldots 2k-1}] .  $$
\end{proof}
By the above Lemma we see that we can write every class which shows up in Lemma \ref{lemma_cap_thom_explicit} as a class in the image of $(p_{V,m})_!$.
Thus, the commutative diagram in Lemma \ref{lemma_comm_diag_gysin} enables us to express the classes in the coproduct through the generators $A_k^i$ and $B_k^i$.
Thus we need to understand the effect of the map
$$  \omega\colon  \mathrm{H}_{\bullet}(V) \xrightarrow[]{\varphi_*} \mathrm{H}_{\bullet}(SM\times_M SM) \to \mathrm{H}_{\bullet}(SM\times SM)        $$
on the classes $[{a}_i]$ and $[{a}_i {b}]$.
Note that these classes can be described as push-forward of the classes $[{a_i}]$ and $[{a_i}{b}]$ in $\mathrm{H}_{\bullet}(SM)$ under the embedding $SM\hookrightarrow V$ given by
$$    [g] \mapsto [g,e] \in V,\quad \text{where}\,\,\, [g]\in G/K_{\gamma} \cong SM ,    $$
see also Lemma \ref{lemma_cohomology_fiber_sm}.
The composition
$$  SM\hookrightarrow V \xrightarrow[]{\varphi} SM\times_M SM \hookrightarrow SM\times SM      $$
is just the diagonal map $d\colon SM\to SM\times SM$.
Hence, we obtain
$$   \omega ([{a}_i]) = d_* [{a_i}] \quad \text{and}\quad \omega ([{a}_i {b}]) = d_*[{a_i}{ b}]      $$
where $d\colon SM\to SM\times SM$ is the diagonal map.
Via the cup ring of $SM$, it is easy to figure out the effect of the diagonal map in homology.
We have
$$   d_*[{a_i}] = \sum_{j=0}^i [{a_j}]\times [{a}_{i-j}]   \quad \text{and}\quad   d_*[{a_i}{ b}] =  \sum_{j=0}^i [{a_j}]\times [{a}_{i-j}{b}] + [{a_j} {b}]\times [{a}_{i-j} ] .       $$

We obtain the following final result.
\begin{theorem} \label{theorem_coproduct}
Let $\mathbb{K}$ be $\mathbb{C}$ or $\mathbb{H}$ and consider $M = \mathbb{K}P^n$.
The string topology coproduct on $M$ behaves as follows.
We have
$$  \cpr A_k^i = \sum_{m=1}^{k-1} \sum_{j=0}^i A_m^j\times A_{k-m}^{i-j}           $$
and
$$    \cpr B_k^i = \sum_{m=1}^{k-1} \sum_{j=0}^i ( B_m^j \times A_{k-m}^{i-j} -A_m^j \times B_{k-m}^{i-j} ) .       $$
\end{theorem}
\begin{proof}
If $Y = A_{k}^i$, we have
$$    Y = (f_k)_* ( - [a_ix_{1\ldots 2k-1}]) .   $$
By Lemmas \ref{lemma_thom_class_k_cycle} and \ref{lemma_cap_thom_explicit}, we have
$$   \tau_{\Lambda}\cap( Y\times [I]) =   \sum_{m=1}^{k-1} (f_k)_*( -[a_i x_{1\ldots 2m-1 \,2m+1\ldots 2k-1}]) \times [\tfrac{2m}{2k}] .    $$
Then by Lemmas \ref{lemma_comm_diag_gysin} and \ref{lemma_gysin_v} we see that
$$   \mathrm{cut}_* \big( \sum_{m=1}^{k-1} (f_k)_*( -[a_i x_{1\ldots 2m-1 \,2m+1\ldots 2k-1}]) \times [\tfrac{2m}{2k}] \big ) = \sum_{m=1}^{k-1} (f_m,f_{k-m})_*\big( (p_{L,m},p_{L,k-m})_!  d_*[{a_i}]  \big) .
$$
Hence we are left with computing the Gysin map of the retraction $(p_{L,m},p_{L,k-m})$.
By \cite[Proposition VI.14.3]{bredon:2013} we see that
$$  (p_{L,m},p_{L,k-m})_! (x\times y) = (-1)^{(\mathrm{dim}(SM) + \mathrm{dim}(\Gamma_m)) (\mathrm{dim}(\Gamma_{k-m})- \mathrm{deg}(y))} \, (p_{L,m})_! (x) \,\times \, (p_{L,k-m})_! (y)       $$
for homology classes $x,y\in\mathrm{H}_{\bullet}(SM)$.
Noting that
$$  \mathrm{dim}(SM)  + \mathrm{dim}(\Gamma_m) = 2(N-1) + \lambda_m    $$
is odd for all $m\in\mathbb{N}$ and that $\mathrm{dim}(\Gamma_{k-m})$ is even for all $m,k\in\mathbb{N}$, $m<k$ we can figure out the sign.
The computation of $\cpr B_k^i$ is analogous.
\end{proof}

\begin{remark}\label{remark_why_figure_eight}
   Note at this point that it is sufficient to consider the critical manifolds with respect to the length function $\mathcal{L}_2$ of the form
    $$  \Sigma_1\times_M \Sigma_{k-1} , \ldots, \Sigma_{k-1}\times_M \Sigma_1 \subseteq \Lambda M\times_M \Lambda M   . $$
    There are two other connected components at level $kl$, namely $M\times_M \Sigma_k$ and $\Sigma_k\times M$.
    However, as we saw now the corresponding homology classes in $\mathrm{H}_{\bullet}(\Lambda M\times_M\Lambda M)$ do not show up in the process of computing the coproduct, therefore we did not consider these components.
\end{remark}


\section{The cohomology product} \label{sec_cohomology}

In this section we briefly describe the Goresky-Hingston product on $\KK P^n$.
We take cohomology with rational coefficients.
Define classes
$$    \sigma_k^i \in \mathrm{H}^{\bullet}(\Lambda M,M) \quad \text{and} \quad \mu_k^i \in \mathrm{H}^{\bullet}(\Lambda M,M)     $$
as follows.
For $k\in \NN$ and $i\in\{0,\ldots,n-1\}$ the class $\sigma_k^i$ is defined to be the dual of $A_k^i$ and the class $\mu_k^i$ is defined to be the dual of $B_k^i$.
Consequently,
$$    \mathrm{deg}(\sigma_k^i) = \lambda_k + i(\lambda +1 ) \quad \text{and} \quad \mathrm{deg}(\mu_k^i) = \lambda_k + (i+1)(\lambda+1) + N-1 .     $$
In particular, note that $\mathrm{deg}(\sigma_k^i)$ is odd and $\mathrm{deg}(\mu_k^i)$ is even for all $k\in\mathbb{N}$, $i\in\{0,\ldots,n-1\}$.
It is clear that the set $$\{\sigma_k^i \,|\, k\in\NN, i\in\{0,\ldots,n-1\}\} \,\, \cup \,\, \{\mu_k^i \,|\, k\in\NN, i\in\{0,\ldots,n-1\}\} $$
generates the cohomology $\mathrm{H}^{\bullet}(\Lambda M,M)$ additively.
Using Theorem \ref{theorem_coproduct} we see that 
$$   \sigma_k^i \ostar \sigma_l^j =   \sigma_{k+l}^{i+j}, \quad \text{if} \,\,\, i+j \leq n-1 \quad \text{and} \quad  \sigma_k^i \ostar \sigma_l^j = 0 \quad \text{else}        $$
and 
$$    \mu_l^j \ostar  \sigma_k^i  =   \mu_{k+l}^{i+j}, \quad \text{if} \,\,\, i+j \leq n-1 \quad \text{and} \quad  \mu_l^j \ostar  \sigma_k^i = 0 \quad \text{else}        $$
for $k,l\in \NN$ and $i,j\in\{0,\ldots,n-1\}$.
Moreover, we have
$$   \mu_k^i\ostar \mu_l^j = 0 \quad \text{for}\,\,\,k,l\in\NN,\,\, i,j\in\{0,\ldots,n-1\} .    $$
We want to mention at this point that the Goresky-Hingston product satisfies the following commutativity property.
If $x\in\mathrm{H}^i(\Lambda M,M)$ and $y\in\mathrm{H}^j(\Lambda M,M)$, then
$$  x\ostar y = (-1)^{(i + \mathrm{dim}(M))(j+\mathrm{dim}(M)) + 1} \, y\ostar x ,  $$
see \cite[Theorem 2.14]{hingston:2017}.
In particular, we get
$$   \sigma_k^i \ostar \sigma_l^j =  \sigma_l^j  \ostar \sigma_k^i \qquad \text{and} \qquad   \sigma_k^i \ostar \mu_l^j  = - \mu_l^j \ostar  \sigma_k^i   $$
which is consistent with the signs in Theorem \ref{theorem_coproduct}.

\begin{theorem} \label{theorem_gh_product}
Let $\mathbb{K}$ be $\mathbb{C}$ or $\mathbb{H}$ and consider the projective space $M = \mathbb{K}P^n$.
The Goresky-Hingston ring $(\mathrm{H}^{\bullet}(\Lambda M, M),\ostar)$ is multiplicatively generated by the classes
$$  \sigma_1^0,\ldots,\sigma_1^{n-1} \quad \text{and} \quad \mu_1^0,\ldots, \mu_1^{n-1}  $$
whose products are subject to the above relations.
In particular, the ring is finitely generated and the element $\sigma_1^0$ is non-nilpotent.
\end{theorem}
\begin{remark}
The existence of a non-nilpotent element is already shown in \cite[Theorem 14.2]{goresky:2009}.
The fact that the Goresky-Hingston ring is finitely generated is analogous to behaviour of the Goresky-Hingston product on spheres, see \cite{goresky:2009}.
\end{remark}


\appendix
\section{Proof of Lemma \ref{lemma_thom_class_k_cycle}}\label{appendix_thom_class}

In this appendix we prove Lemma \ref{lemma_thom_class_k_cycle}.
We will use the notation established throughout the paper.
We shall show that the pullback of the class $\tau_{\Lambda}\in \mathrm{H}^{N}(U_{\Lambda},U_{\Lambda,\geq\epsilon_0})$ under the map 
$$   (f_k,\id_I)\colon (U_{\Gamma_k},U_{\Gamma_k,\geq\epsilon_0}) \to (U_{\Lambda},U_{\Lambda,\geq\epsilon_0})     $$
is given by 
$$  \tau_k =  (f_k,\id_I)^*\tau_\Lambda = \sum_{m=1}^{k-1}   \xi_{2m}\times \eta_m .      $$

First, let us consider the generators of $\mathrm{H}^{N-1}(\Gamma_k)$.
The classes
$$   \xi_2,\ldots,\xi_{2k-2}        $$
are generators of $\mathrm{H}^{N-1}(\Gamma_k)$, so are the classes 
$$    \alpha^m\cup \xi_{i_1} \cup \ldots \cup \xi_{i_{l}} \in \mathrm{H}^{N-1}(\Gamma_k)\quad \text{with}\quad m\in\{0,\ldots,n-1\},\,\,\, 1\leq i_1< \ldots i_{l} \leq 2k-1, \,\, i_j \,\,\text{odd}       $$
where necessarily 
$$   m (\lambda + 1) + l \, \lambda = N-1 .   $$
One can check that these are the only generators.
Hence, we have
$$    \tau_k = \sum_{j=1}^{k-1}\sum_{m=1}^{k-1} \lambda_{j,m} \xi_{2j}\times \eta_m + \sum_{m, \,\,1\leq i_1<\ldots < i_l \leq 2k-1 }\sum_{s=1}^{k-1} \rho_{m,i_1\ldots i_l,s}( \alpha^m\cup \xi_{i_1} \cup \ldots \cup \xi_{i_{l}} )\times \eta_s        $$
where $\lambda_{j,m}\in\mathbb{Q}$ and $\rho_{m,i_1,\ldots,i_l,s}\in\mathbb{Q}$ are coefficients and where the sum in the second term is taken over those combinations of $m,i_1,\ldots,i_l$ such that
$$   m(\lambda + 1) + l\,\lambda = N-1 \quad \text{and all} \,\,i_j \,\,\text{odd}.   $$

We want to show that 
$$  \lambda_{j,m} = 1 ,\quad \text{if}\,\,\, j=m\in\{1,\ldots,k-1\}   \quad \text{and}\quad \lambda_{j,m} = 0 \,\,\,\text{else}   $$
and that all coefficients $\rho_{m,i_1\ldots i_l,s} = 0$.
\begin{claim}
    The coefficients $\lambda_{j,m}\in\mathbb{Q}$ satisfy
    $$  \lambda_{j,m} = 1 ,\quad \text{if}\,\,\, j=m  \quad \text{and}\quad \lambda_{j,m} = 0 \,\,\,\text{else}   . $$
\end{claim}
\begin{proof}
Fix $j,m\in\{1,\ldots,k-1\}$.
We have
\begin{equation}\label{eq_lambda_jm}
        \lambda_{j,m} = \langle  \tau_k, [x_{2j}]\times [I_m] \rangle .    
\end{equation}
Recall that $\tau_{\Lambda} = \mathrm{ev}_{\Lambda}^*\tau_M$.
Hence, we get
\begin{equation}\label{eq_lambda_jm_2}
     \langle \tau_k,[x_{2j}]\times [I_m]\rangle = \langle (f_k,\id_I)^*\mathrm{ev}_{\Lambda}^*\tau_M , [x_{2j}]\times [I_m]\rangle = \langle \tau_M , (\mathrm{ev}_{\Lambda}\circ (f_k,\id_I))_* ([x_{2j}]\times [I_m])\rangle   .    
\end{equation}
Thus we consider
$$   g_m = \mathrm{ev}_{\Lambda} \circ (f_k,\id_{I_m}) : \Gamma_k \times (I_m,J_m) \to (U_M,U_{M,\geq\epsilon_0})   .   $$
Recall that the class $[x_{2j}]$ can be described as follows.
Let $s_{2j}\colon K/K_{\gamma}\hookrightarrow \Gamma_k$ be the embedding as constructed in Lemma \ref{lemma_embed_sub_k_cycle}, i.e.
$$   s_{2j}([x]) = [e,\ldots, e,x,e\ldots, e],\quad \text{for} \,\,\, [x]\in K/K_{\gamma} .  $$
Here, the $x$ appears at position $2j$.
Then $[x_{2j}] = (s_{2j})_*[\mathbb{S}^{N-1}] $, where $[\mathbb{S}^{N-1}]$ is the fundamental class of $K/K_{\gamma} \cong \mathbb{S}^{N-1}$.
Let $[C]\in \mathrm{H}_{N-1}(U_M,U_{M,\geq\epsilon_0} )$ be the class dual to $\tau_M$. We want to show that
$$   (g_m)_*([x_{2j}]\times [I_m]) =  [C]  \,\,\,\,\text{if and only if} \,\,\, j = m \,\,\,\,\text{and } 0 \,\,\text{otherwise} .   $$ 
Set 
$$  B_{p_0} = \{q\in M\,|\,\mathrm{d}(p_0,q)<\epsilon\} \quad \text{and}\quad B_{p_0,\geq\epsilon_0} = \{q\in B_{p_0}\,|\,\mathrm{d}(q,p_0)\geq\epsilon_0\} .    $$
Note that if 
$$  i\colon (B_{p_0},B_{p_0,\geq\epsilon_0}) \hookrightarrow (U_M, U_{M,\geq\epsilon_0}) ,\qquad  i(q) = (p_0,q)    $$
for $q\in B_{p_0}$
is the inclusion of the fiber of the normal tubular neighborhood then we have
$$   [C] = i_* [B], \quad \text{where}\quad [B]\in \mathrm{H}_{N-1}(B_{p_0},B_{p_0,\geq\epsilon_0})    $$
is a positively oriented generator.
We define the map 
$$h_m = g_m \circ (s_{2j},\id_{I_m})  \colon  K/K_{\gamma} \times (I_m,J_m)\to (U_M,U_{M,\geq \epsilon_0}) .  $$
Note that this factors through maps
$$ \begin{tikzcd}
    (K/K_{\gamma}\times I_m, K/K_{\gamma}\times J_m)  \arrow[swap]{rd}{h_m'} \arrow[]{rr}{h_m} & & ( U_M,U_{M,\geq\epsilon_0})
    \\
     & (B_{p_0},B_{p_0,\geq\epsilon_0}) \arrow[swap]{ru}{i} 
\end{tikzcd}   $$
Define $\mathrm{ev}\colon \Lambda\times I\to M$ by $\mathrm{ev}(\gamma,s) = \gamma(s)$.
In order to show that
$$ (g_m)_*([x_{2m}]\times [I_m]) =  [C]    $$
it thus suffices to shows that
$$   (h_m')_*([x_{2m}]\times [I_m]) =  [B]   $$
where $(g_m)' = \mathrm{ev}\circ (f_k,\id_{I_m})  $.
We have that $h_m' = g_m'\circ (s_{2j},\id_{I_m})$.
Recall that the orbit of a point $q\in B_{p_0}$ under $K$ is the distance-sphere around $p_0$ of radius $\mathrm{d}(p_0,q)$. 
We compute the map $h_m'$ explicitly.
Let $x\in K$.
In case that $m < j$, we have
\begin{equation} \label{eq_evaluation_1}
     h_m' ([x],t)    =\begin{cases} 
      \gamma(t) & t < \tfrac{2m}{2k}  \\
      \gamma(t)  & t\geq  \tfrac{2m}{2k}
   \end{cases}
\end{equation}  
in case $m = j$, we get
\[ h_m' ([x],t)    =\begin{cases} 
     \gamma(t) & t < \tfrac{2m}{2k}  \\
      x.\gamma(t)  & t\geq \tfrac{2m}{2k}
   \end{cases}
\]  
and in case $m > j$, we obtain
\[ h_m' ([x],t)    =\begin{cases} 
      x.\gamma(t) & t < \tfrac{2m}{2k}  \\
      x.\gamma(t)  & t\geq \tfrac{2m}{2k} .
   \end{cases} 
\]  
Now, consider the following commutative diagram
$$ 
\begin{tikzcd}
    \mathrm{H}_{N-1}(\mathbb{S}^{N-1} )\otimes \mathrm{H}_1(I_m,J_m)  \arrow[]{r}{\id\otimes\partial} \arrow[]{d}{\cong}
    & [2.5em] 
    \mathrm{H}_{N-1}(\mathbb{S}^{N-1})\otimes \mathrm{H}_0(J_m)  \arrow[]{d}{\cong}
    \\ 
    \mathrm{H}_{N}( \mathbb{S}^{N-1}\times I_m,\mathbb{S}^{N-1}\times J_m   ) \arrow[]{r}{\partial}  \arrow[]{d}{(h'_m)_*} 
    &
    \mathrm{H}_{N-1}(\mathbb{S}^{N-1}\times J_m) \arrow[]{d}{(h_m')_*}
    \\
    \mathrm{H}_{N}(B_{p_0},B_{p_0,\geq\epsilon_0}) \arrow[]{r}{\cong} \arrow[]{d}{\cong} 
    & 
    \mathrm{H}_{N-1}(B_{p_0,\geq\epsilon_0}) \arrow[]{d}{\cong}
    \\
    \mathrm{H}_{N}(\mathbb{D}^{N},\mathbb{S}^{N-1}) \arrow[r, "\cong"', "\partial"]
    &
    \mathrm{H}_{N-1}(\mathbb{S}^{N-1}) 
\end{tikzcd}
$$
where the maps $\partial$ are the respective connecting homomoprhisms.
The middle square is
$$   
\begin{tikzcd}
    \QQ  \arrow[]{r}{  x \mapsto (-x, x)  } \arrow[swap]{d}{ h_{m,j}^1   }
    &  [3em]
    \QQ \oplus \QQ \arrow[]{d}{ h_{m,j}^2 }
    \\ 
    \QQ \arrow[]{r}{\id} 
    &
    \QQ 
\end{tikzcd}
$$  
with $h_{m,j}^1$ and $h_{m,j}^2$ the maps induced by $(h_m')_*$.
In case $m < j$, it is clear that $h_{m,j}^2$ is the trivial map, so $h_{m,j}^1 = 0$.
This can be seen from equation \eqref{eq_evaluation_1} since the map $h_m'$ is homotopic to a locally constant map.
If $m > j$, we have by the orientation convention of Proposition \ref{prop_cup_ring_gamma_k} that $h_{m,j}^2(x,y) = x+ y$, so again $h_{m,j}^1$ is trivial.
Finally, for $m = j$, we see that
$h_{m,m}^2 (x,y) = y$, so we get
$$   h_{m,m}^1 = \id  . $$
Here we use again the orientation convention established in Proposition \ref{prop_cup_ring_gamma_k}.
This shows that $$(h_m')_* ([x_{2m}] \times [I_m]) = [B]\quad \text{and}\quad (h_m')_*([x_{2j}]\times [I_m]) = 0 \quad  \text{for} \,\, j\neq m. $$
The claim then follows from equations \eqref{eq_lambda_jm} and \eqref{eq_lambda_jm_2}.
\end{proof}


\begin{claim}
    All coefficients $\rho_{m,i_1\ldots i_l,s}$ vanish.
\end{claim}
\begin{proof}
Fix $m\in\{0,\ldots,n-1\}$ and odd integers $i_1,\ldots,i_l$ with $1\leq i_1< \ldots <i_l\leq 2k-1$ such that
$$    m\,(\lambda + 1) +  l\,\lambda = N-1 .       $$
We also fix $s\in\{1,\ldots,k-1\}$.
We begin by describing a dual class to the cohomology class
$$  \alpha^m \cup \xi_{i_1}\cup \ldots\cup \xi_{i_l} \in\mathrm{H}^{N-1}(\Gamma_k) .   $$

Let $\kappa = m(\lambda + 1) $.
Recall that $\pi_L \colon L_k\cong SM\to M$ is the unit sphere bundle of the underlying manifold.
Since $\kappa < N$ we see from the Gysin sequence of this sphere bundle that there is an isomorphism
$$  (\pi_{L})_* \colon \mathrm{H}_{\kappa}(SM) \xrightarrow[]{\cong} \mathrm{H}_{\kappa}(M) .     $$
We know that $\mathrm{H}_{\kappa}(M)\cong\mathbb{Q}$.
Note that a generator of $\mathrm{H}_{\kappa}(M)$ can be described as follows.
It is well-known that there is an inclusion of $\mathbb{K}P^m$ into $M = \mathbb{K}P^n$ which maps a fundamental class of $\mathbb{K}P^m$ to a generator of $\mathrm{H}_{\kappa}(M)$.
In particular we can choose this inclusion in such a way that the basepoint $p_0\in M$ is also the basepoint of $\mathbb{K}P^m$.
Denote this inclusion by $j\colon \mathbb{K}P^m \to \mathbb{K}P^n$.
We pull back the unit sphere bundle along $j$.
Since the dimension of the fiber of this bundle is greater than the dimension of the base the Euler class vanishes and hence this bundle has a section $s\colon \mathbb{K}P^m\to j^*SM$.
We compose this with the canonical map $j^*SM \to SM$ to get a section $\sigma\colon\mathbb{K}P^m\to SM$.
It is clear that we have
$$   \pi_L \circ \sigma = j,    $$
so we see that we obtain an isomorphism
$$   \sigma_* \colon \mathrm{H}_{\kappa}(\mathbb{K}P^m) \xrightarrow[]{\cong} \mathrm{H}_{\kappa}(SM) .   $$
Hence, we can represent a generator in degree $\kappa$ by the image of a fundamental class of $\mathbb{K}P^m$ under the map $\sigma$.

Consider the manifold $\Gamma^{0i_1\ldots i_l}$ where we use the notation of Lemma \ref{lemma_embed_sub_k_cycle}.
We have
$$  \Gamma^{0i_1\ldots i_l} =  G\times_{K_{\gamma}} K_a \times_{K_{\gamma}} \ldots \times_{K_{\gamma}} (K_a/K_{\gamma}) .   $$
In Lemma \ref{lemma_embed_sub_k_cycle} we saw that this embeds into the manifold $\Gamma_k$ via a map $s_{0i_1\ldots i_l}\colon \Gamma^{0i_1\ldots i_l}\to \Gamma_k$.
Moreover, as in the discussion after the proof of Lemma \ref{lemma_embed_sub_k_cycle} one sees that this manifold is an iterated sphere bundle over $SM$ and the fiber at each step of the iterated sphere bundle is the sphere $\mathbb{S}^{\lambda}$.
We now pull back this fiber bundle along the embedding $\sigma\colon \mathbb{K}P^m\hookrightarrow SM$ to get a pull-back space $X$ which is clearly a manifold.
We have a commutative diagram
$$     
\begin{tikzcd}
    X \arrow[hook]{r}{} \arrow[]{d}{} \arrow[rr, bend left, "\iota_m"]
    &
    \Gamma^{0i_1\ldots i_l } \arrow[hook, swap]{r}{s_{0i_1\ldots i_l}} \arrow[]{d}{}
    &
    \Gamma_k \arrow[]{d}{}
    \\
    \mathbb{K}P^m \arrow[]{r}{\sigma} 
    &
    SM \arrow[]{r}{=}
    & SM
\end{tikzcd}
$$
where the left square is just the pullback diagram.
The manifold $X$ is also an iterated sphere bundle with base $\mathbb{K}P^m$.
Moreover, from the Gysin sequences of $X$, $\Gamma^{0i_1\ldots i_l}$ and $\Gamma_k$ one can see that $(\iota_m)_* [X]$ is indeed the dual homology class to the cohomology class $\alpha^m\cup \xi_{i_1}\cup\ldots \cup \xi_{i_l}$.
Consequently, we obtain 
\begin{eqnarray} \label{eq_rho_coefficients}
          \rho_{m,i_1\ldots i_l,s}  & =  & \langle \tau_k ,  (\iota_m,\id_{I_s})_* ([X]   \times [I_s] )   \rangle    \\
          &=& \nonumber
          \langle \tau_M,( g\circ (\iota_M,\id_{I_s}))_* ([X]\times [I_s]) \rangle 
\end{eqnarray}
where we have $g = \mathrm{ev}_{\Lambda}\circ (f_k,\id_{I_s})$ as before.
We define $$h_{m,s} = g\circ (\iota_M,\id_{I_s}) \colon X\times (I_s,J_s)\to (U_M, U_{M,\geq\epsilon_0}) . $$
If we show that the class
$$   ( h_{m,s})_* ([X]\times [I_s]) \in \mathrm{H}_{N}(U_M,U_{M,\geq \epsilon_0})   $$
vanishes, then by equation \eqref{eq_rho_coefficients} we have shown that the coefficients $\rho_{m,i_1,\ldots,i_l,s}$ vanish.
Let us compute the map $h_{m,s}$ explicitly.
There is a number $o \in \{ 1,\ldots,l\}$ such that 
$$   i_1< \ldots < i_{o} < 2s < i_{o +1} < \ldots i_l  .   $$
Let $[g,k_1,\ldots,k_l]\in \Gamma^{0i_1\ldots i_l}$ with $\sigma(x) = [g]\in G/K_{\gamma}\cong SM$ for some $x\in \mathbb{K}P^m$.
Furthermore, let $t\in I_s$, then we have
$$    h_{m,s} ([g,k_1,\ldots,k_l],t) =  ( g.\gamma(0) ,  gk_1\ldots k_o . \gamma(t) ) = (x, gk_1\ldots k_o . \gamma(t) )  .     $$
We see from the above expression that the map $h_{m,s}$ factors as follows
$$ 
\begin{tikzcd}
    X\times (I_s,J_s) \arrow[]{rr}{h_{m,s}} \arrow[swap]{rd}{h_{m,s}'}
    & & 
    (U_M,U_{M,\geq\epsilon_0})
    \\
    &
    (t^* U_M, t^* U_{M,\geq\epsilon_0})  \arrow[hook]{ru}{}   &    
\end{tikzcd}
$$
where $t\colon \mathbb{K}P^m\hookrightarrow \mathbb{K}P^n = M$ is the inclusion.
We now consider the following commutative diagram
$$
\begin{tikzcd}
&
\mathrm{H}_{N-1}(X)\otimes \mathrm{H}_1(I_s,J_s) \arrow[]{r}{\id\otimes \, \partial} \arrow[]{d}{\cong}
& [2.em]
\mathrm{H}_{N-1}(X)\otimes \mathrm{H}_0(J_s) \arrow[]{d}{\cong}.
&
\\
    &
    \mathrm{H}_N(X\times I_s, X\times J_s) \arrow[]{r}{\partial} \arrow[]{d}{(h_{m,s}')_*}
    &
    \mathrm{H}_{N-1}(X\times J_s) \arrow[]{d}{(h_{m,s}')_*}
    &
    \\
    \mathrm{H}_{N}(t^* U_M) \arrow[]{r}{}
    &
    \mathrm{H}_N(t^* U_M,t^* U_{M,\geq\epsilon_0}) \arrow[]{r}{\partial}
    &
    \mathrm{H}_{N-1}(t^* U_{M,\geq\epsilon_0}) \arrow[]{r}{} 
    & \mathrm{H}_{N-1}(t^* U_M) .
\end{tikzcd}
$$
Note that the space $t^* U_M$ is homeomorphic to a disk bundle over $\mathbb{K}P^m$ and is therefore homotopy equivalent to $\mathbb{K}P^m$ itself.
Clearly, the dimension of $\mathbb{K}P^m$ satisfies $\mathrm{dim}(\mathbb{K}P^m) \leq N-2$ and therefore the homology groups on the very left and the very right in the lower row of the above diagram vanish.
Therefore the connecting homomorphism in the lower row is an isomorphism.

Moreover, note that $J_s$ is adisjoint union of two intervals, hence it is homotopy equivalent to a union of two points, i.e. we have a homotopy equivalence
$$  \Theta\colon    X\times \{ t_-\} \cup X\times \{t_+\} \xrightarrow[]{\hphantom{b}\simeq\hphantom{b}} X \times J_s    $$
where $t_- \in (\tfrac{2s}{2k} - \delta, \tfrac{2s}{2k}  + \delta_0] $ and $t_+ \in [\tfrac{2s}{2k} + \delta_0, \tfrac{2s}{2k} + \delta)$.
We can choose $t_-$ and $t_+$ such that they are equidistant to $\tfrac{2s}{2k}$, i.e. $|t_+ - \tfrac{2s}{2k}| =|t_- - \tfrac{2s}{2k}| $.
This implies that 
$$ \delta_1 := \mathrm{d}(\gamma(t_+),\gamma(0)) = \mathrm{d}(\gamma(t_-),\gamma(0)) .  $$


We thus consider the maps
$$  k_{m,s,-} \colon X\times \{t_-\} \xrightarrow[]{} X\times J_s \xrightarrow[]{h'_{m,s}} t^*U_{M,\geq\epsilon_0}  
$$
and 
$$  k_{m,s,+} \colon X\times \{t_+\} \xrightarrow[]{} X\times J_s \xrightarrow[]{h'_{m,s}} t^*U_{M,\geq\epsilon_0}  
. $$
Note that $X\cong X\times \{t_{\pm}\}$, so we can understand both maps as maps $X\to 
t^* U_{M,\geq\epsilon_0} $.
We now show that they induce the same map in homology.

First, we define a map $\varphi\colon X\to X$.
Note that since $M$ is a symmetric space there is an isometry $S\colon M\to M$ which fixes the basepoint $p_0 = \gamma(0)$ and acts as $-\mathrm{id}_{T_{p_0} M}$ on the tangent space.
This isometry reverses geodesics going through the basepoint.
Since the point $a\in M$ is the unique conjugate point in the interior of the prime closed geodesic $\sigma$ we have
$$   S.a = S.\gamma(\tfrac{1}{2}) = \gamma(-\tfrac{1}{2}) = a .   $$
Consequently, the isometry satisfies $S\in K_a$.
Therefore, we can define a map
$$  \Phi\colon G\times (K_a)^l \to G\times (K_a)^l     $$
by setting
$$    \Phi(g,k_1,\ldots,k_l) = (g,k_1,\ldots, k_l S) \quad \text{for} \quad g\in G, \,k_1,\ldots,k_l\in K_a .    $$
Note that the element $S$ commutes with all elements in $K$.
This is because the isotropy representation of a symmetric space is faithful and the element $-\mathrm{id}\in O(N)$ is clearly in the center of $O(N)$.
Hence, one sees that the map $\Phi$ is indeed equivariant with respect to the $(K_{\gamma})^{l+1}$-action and therefore induces a smooth map $\varphi' \colon \Gamma^{0i_1\ldots i_l}\to \Gamma^{0i_1\ldots i_l}$.
Since this map squares to the identity it is a diffeomorphism.
Moreover it clearly respects the fiber bundle structure, so it restricts to a map $\varphi\colon X\to X$.
We consider this map for the following reason.
Since the isometry $S$ reverses geodesics through the basepoint we have
$$      S. \gamma(t_-) = \gamma(t_+)            $$
by our choice of $t_-$ and $t_+$.
Therefore we see by definition of $k_{m,s,-}$ and $k_{m,s,+}$ that
$$   k_{m,s,+} = k_{m,s,-}\circ \varphi .     $$
Hence, if we show that $\varphi_* \colon \mathrm{H}_{\kappa}(X) \to \mathrm{H}_{\kappa}(X)$ is the identity, it follows that $k_{m,s,+}$ and $k_{m,s,-}$ induce the same map in homology.
We shall argue that the degree of $\varphi$ is $1$.
Since $\varphi$ is a diffeomorphism it suffices to check whether the differential at a given point is orientation-preserving or orientation-reversing.
Let $[e,e,\ldots,e]\in \Gamma^{0i_1\ldots i_l}$.
There is an open neighborhood of $[e,\ldots,e]\in \Gamma^{0i_1\ldots i_{l-1}}$ such that
the fiber bundle $p\colon \Gamma^{0i_1\ldots i_l}\to \Gamma^{0i_1\ldots i_{l-1}}$ is trivial over $U$, i.e.
$$     p^{-1}(U) \cong U \times \mathbb{S}^{\lambda} .  $$
But in this local trivialization it is very easy to understand the effect of the map $\varphi$.
We have
$$    \varphi|_{p^{-1}(U)} \colon p^{-1}(U) \to p^{-1}(U) \quad \text{satisfies} \quad \varphi(u,x) = (u,-x)     $$
for $u\in U$, $x\in \mathbb{S}^{\lambda}$.
Now, the identity on $U$ is clearly orientation-preserving as is the antipodal map on an odd-dimensional sphere.
Therefore we get that $\varphi_* = \id_{\mathrm{H}_{\kappa}(X)}$.

Now, consider again the class
$$  [X]\times [I_s] \in \mathrm{H}_N(X\times I_s,X\times J_s) .    $$
It is well-known that the connecting homomorphism maps this to
$$  \partial ([X]\times [I_s]) =  [X]\times [t_+] - X\times [t_-] .  $$
But as we have seen now
$$   (h'_{m,s})_* ([X] \times [t_+]) =   (k_{m,s,+})_*([X]) =     (k_{m,s,-})_*([X])  =  (h'_{m,s})_* ([X] \times [t_-])  $$
so this shows that
$$ (h'_{m,s})_* \circ \partial \, ([X]\times [I_s]) = 0 .$$
Consequently, we obtain $\rho_{m,i_1\ldots i_l,s} = 0$.
\end{proof}
Note that the strategy of the proof of the last claim is very similar to the methods employed in \cite[Section 7]{kupper:2022}.
The proof of the two claims completes the proof of Lemma \ref{lemma_thom_class_k_cycle}.

\section{Relative Cap Product} \label{appendix_cap}

In this section we review the construction of the relative cap product which is used in the definition of the string topology coproduct.
We closely follow \cite[Appendix A]{hingston:2017}.

Assume that $X$ is a topological space with subspaces $A,B\subseteq X$ such that
$$  \mathrm{C}_{\bullet}(A) + \mathrm{C}_{\bullet}(B) \hookrightarrow \mathrm{C}_{\bullet}(A\cup B)     $$
is a quasi-isomorphism.
Then there is a relative cap-product
\begin{equation} \label{eq_cap_product}
 \cap : \mathrm{H}^k(X,A)\otimes \mathrm{H}_m(X,A\cup B) \to \mathrm{H}_{m-k}(X,B)  .   \end{equation} 
The condition on the subspaces is satisfied if e.g. both $A$ and $B$ are open. See \cite[Section VI.5]{bredon:2013} for details.

Assume now that $U_0\subseteq U_1\subseteq X$ are open subsets of $X$ such that $\mathcal{U} = \{U_1, \mathrm{int}(U_0^c)\}$ is an open cover of $X$.
Here we use the notation $U_0^c = X\setminus U_0$.
Furthermore, assume that $A\subseteq X$ is another subset which is not necessarily required to be a subset of $U_0$ or $U_1$.
We assume that the intersections $U_1\cap U_0^c$ and $U_1\cap A$ are such that the relative cap-product 
$$  \cap\colon \mathrm{H}^k(U_1,U_1\cap U_0^c ) \otimes \mathrm{H}_mk(U_1,U_1\cap U_0^c \cup U_1\cap A) \to \mathrm{H}_{m-k}(U_1,U_1\cap A)    $$
as in equation \eqref{eq_cap_product} is defined.
Then if $t\in \mathrm{C}^k(U_1,U_1\cap U_0^c)$,
we define a map
$$   t \, \cap'  : \mathrm{C}_m(X,A) \to \mathrm{C}_{m-k}(U_1,U_1\cap A)      $$
as the composition of the maps
\begin{eqnarray*} 
      \mathrm{C}_{m}(X,A) \xrightarrow{} \mathrm{C}_{m}(X, U_0^c \cup A)
     \xrightarrow{ \rho } \mathrm{C}_{m}^{\mathcal{U}}( X,U_0^c\cup A )  \xrightarrow{} \mathrm{C}_{ m }(  U_1, U_1\cap U_0^c \cup U_1\cap A  ) \\ 
     \xrightarrow{ t\, \cap \,\, } 
    \mathrm{C}_{m-k}( U_1, U_1\cap A)
\end{eqnarray*}
where $\rho$ is a map that subdivides chains with respect to the open cover $\mathcal{U}$, e.g. barycentric subdivision, see e.g. \cite[Section IV.17]{bredon:2013}.
The cap product in the last step is then a chain-level version of the cap product in equation \eqref{eq_cap_product}.
This composition is a chain map and therefore induces a map in homology
$$  \cap' : \mathrm{H}^k(U_1,U_1\cap U_0^c) \otimes \mathrm{H}_m(X,A) \to \mathrm{H}_{m-k}(U_1,U_1\cap A)   $$
which we will refer to as \textit{cap product} as well.
From now on, we will also denote it by $\cap$ and from the context it will be clear whether we are referring to this cap product or to the one of equation \eqref{eq_cap_product}.
We now state a naturality statement for the cap product.
\begin{prop} \label{prop_nat_cap}
Let $(X,U_1,U_0)$ and $(Y,V_1,V_0)$ be triples of spaces and let $A\subseteq X$ and $B\subseteq Y$ be subsets such that the inclusion $$\mathrm{C}_{\bullet}(U_1\cap U_0^c) + \mathrm{C}_{\bullet}(U_1\cap A) \hookrightarrow \mathrm{C}_{\bullet}(U_1\cap (U_0^c\cup A))$$ is a quasi-isomorphism and similarly for $V_1,V_0$ and $B$.
Then if $f\colon (X,U_1,U_0)\to (Y,V_1,V_0)$ is a map of triples such that $f(A)\subseteq B$ and if $\tau\in \mathrm{H}^k(V_1,V_1\cap V_0^c)$ is a cohomology class then the diagram
$$  
\begin{tikzcd}
\mathrm{H}_m(X,A) \arrow[]{r}{f_*} \arrow[]{d}{f^*(\tau) \cap } & \mathrm{H}_m(Y,B) \arrow[]{d}{\tau\cap} 
\\
\mathrm{H}_{m-k}(U_1,U_1\cap A) \arrow[]{r}{f_*} & \mathrm{H}_{m-k}(V_1,V_1\cap B)
\end{tikzcd}
$$
commutes.
\end{prop}
\begin{proof}
The proof can be done analogously to the one of \cite[Lemma A.1]{hingston:2017} since the restriction $A\subseteq U_0$ in the proof of \cite[Lemma A.1]{hingston:2017} is not necessary.
\end{proof}

We also need a statement about the compatibility of the relative cap product and the usual cross products.
\begin{prop} \label{prop_compat_cross_cap}
Let $X$ and $Y$ be topological spaces.
Furthermore, let $A,U_0,U_1\subseteq X$ be subspaces such that $(X,U_1,U_0)$ is a triple of spaces and such that the inclusion $$\mathrm{C}_{\bullet}(U_1\cap U_0^c) + \mathrm{C}_{\bullet}(U_1\cap A) \hookrightarrow \mathrm{C}_{\bullet}(U_1\cap (U_0^c\cup A))$$ is a quasi-isomorphism.
Let $\xi \in \mathrm{H}^i(Y)$, $\eta\in \mathrm{H}^j(U_1,U_1\cap U^c_0)$, $y\in \mathrm{H}_m(Y)$ and $z\in\mathrm{H}_n(X,A)$.
Then
$$   (\xi\times \eta) \cap (y\times z) = (-1)^{jm} (\xi\cap y)\times (\eta\cap z) .     $$
\end{prop}
\begin{proof}
If $$\mathrm{C}_{\bullet}(U_1\cap U_0^c) + \mathrm{C}_{\bullet}(U_1\cap A) \hookrightarrow \mathrm{C}_{\bullet}(U_1\cap (U_0^c\cup A))$$ is a quasi-isomorphism then clearly the same property holds for
$$\mathrm{C}_{\bullet}(Y\times (U_1\cap U_0^c)) + \mathrm{C}_{\bullet}(Y\times (U_1\cap A)) \hookrightarrow \mathrm{C}_{\bullet}(Y\times (U_1\cap (U_0^c\cup A))) .$$
Consequently, the proof of the analogous property for the usual cap product carries over, see \cite[Theorem VI.5.4]{bredon:2013}.
\end{proof}
\begin{example} \label{ex_relative_cap_interval}
Take $X = I$ to be the unit interval with $A = \partial I = \{0,1\}$.
Furthermore, choose a small $\delta>0$ and a number $\delta_0>0$ with $\delta_0<\delta$.
Clearly, the homology group $\mathrm{H}_1(I,\partial I)$ is generated by a class $[I]$ which is represented by the relative cycle
$$  \sigma \colon I\to I,\quad \sigma = \id_I .      $$
Now, choose $$U_1 = (\tfrac{1}{2}-\delta,\tfrac{1}{2}+\delta )\quad \text{and}\quad U_0 = (\tfrac{1}{2}-\delta_0,\tfrac{1}{2}+\delta_0)  .$$
We have $U^c_0 \cap A = U_0^c$.
Note that the subdivision with respect to the cover $\mathcal{U} = \{ U_1, \mathrm{int}(U_0^c)\}$, i.e. the map
$$  \mathrm{C}_1(I,U_0^c) \to \mathrm{C}_1^{\mathcal{U}}(I,U_0^c)      $$
can be chosen as follows.
It maps $\sigma$ to $\sigma_1+\sigma_2+\sigma_3$
where $$  \sigma_1\colon I\to [0,a_1],\quad \sigma_2\colon I\to [a_1,a_2] \quad \text{and}\quad \sigma_3\colon I\to [a_2,1]    $$
are the respective affine linear maps and where $$a_1\in (\tfrac{1}{2}-\delta,\tfrac{1}{2}-\delta_0) \quad \text{and}\quad  a_2\in (\tfrac{1}{2}+\delta_0, \tfrac{1}{2}+\delta) .     $$
We now want to determine the cap product with a representative of a generator of $$\mathrm{H}^1(U_1,U_1\cap U_0^c) = \mathrm{H}^1\big((\tfrac{1}{2}-\delta,\tfrac{1}{2}+\delta) , (\tfrac{1}{2}-\delta,\tfrac{1}{2}-\delta_0]\cup[\tfrac{1}{2}+\delta_0,\tfrac{1}{2}-\delta)\big) .$$
If we choose a cocycle $\tau\in \mathrm{C}^1(U_1,U_1\cap U_0^c)$ representing a generator of $\mathrm{H}^1(U_1,U_1\cap U_0^c)$ which is dual to $\sigma_2\in\mathrm{C}_1(U_1,U_1\cap U_0^c)$ then it is clear that we get
$$  \tau\cap \sigma_2 = t_0 \in \mathrm{C}_0(U_1)     $$
for some point $t_0\in U_1$.
Therefore in homology, we see that the relative cap product $[\tau]\cap [I]$ is
$$   [\tau]\cap [I] = [t_0]\in \mathrm{H}_1(U_1) .    $$
We use this example in Section \ref{sec_computation}.
\end{example}



\bibliography{lit}

\providecommand{\bysame}{\leavevmode\hbox to3em{\hrulefill}\thinspace}
\providecommand{\MR}{\relax\ifhmode\unskip\space\fi MR }
\providecommand{\MRhref}[2]{%
  \href{http://www.ams.org/mathscinet-getitem?mr=#1}{#2}
}
\providecommand{\href}[2]{#2}
\begin{thebibliography}{NRW22}

\bibitem[Ara62]{araki:1962}
Sh{\^o}r{\^o} Araki, \emph{On bott-samelson $ k $-cycles associated with symmetric spaces}, Journal of Mathematics, Osaka City University \textbf{13} (1962), no.~2, 87--133.

\bibitem[Bes78]{besse:1978}
Arthur~L Besse, \emph{Manifolds all of whose geodesics are closed}, vol.~93, Springer Science \& Business Media, 1978.

\bibitem[Bot56]{bott:1956}
Raoul Bott, \emph{An application of the {M}orse theory to the topology of {L}ie groups}, Bulletin de la Soci{\'e}t{\'e} Math{\'e}matique de France \textbf{84} (1956), 251--281.

\bibitem[Bot58]{bott:1958b}
\bysame, \emph{The space of loops on a lie group.}, Michigan Mathematical Journal \textbf{5} (1958), no.~1, 35--61.

\bibitem[Bre13]{bredon:2013}
Glen~E Bredon, \emph{Topology and geometry}, Graduate Texts in Mathematics, vol. 139, Springer Science \& Business Media, 2013.

\bibitem[BS58]{bott:1958a}
Raoul Bott and Hans Samelson, \emph{Applications of the theory of {M}orse to symmetric spaces}, American Journal of Mathematics \textbf{80} (1958), no.~4, 964--1029.

\bibitem[GH09]{goresky:2009}
Mark Goresky and Nancy Hingston, \emph{Loop products and closed geodesics}, Duke Math. J. \textbf{150} (2009), no.~1, 117--209.

\bibitem[HO13]{hingston:2013oancea}
Nancy Hingston and Alexandru Oancea, \emph{The space of paths in complex projective space with real boundary conditions}, arXiv preprint arXiv:1311.7292 (2013).

\bibitem[HW17]{hingston:2017}
Nancy Hingston and Nathalie Wahl, \emph{Product and coproduct in string topology, revised version 2021}, arXiv preprint arXiv:1709.06839 (2017).

\bibitem[HW19]{hingston:2019}
\bysame, \emph{Homotopy invariance of the string topology coproduct}, arXiv preprint arXiv:1908.03857 (2019).

\bibitem[Kli78]{klingenberg:78}
Wilhelm Klingenberg, \emph{Lectures on closed geodesics}, Grundlehren der Mathematischen Wissenschaften, Vol. 230, Springer-Verlag, Berlin-New York, 1978.

\bibitem[Kli95]{klingenberg:1995}
\bysame, \emph{Riemannian geometry}, second ed., De Gruyter Studies in Mathematics, vol.~1, Walter de Gruyter \& Co., Berlin, 1995.

\bibitem[KS22]{kupper:2022}
Philippe Kupper and Maximilian Stegemeyer, \emph{On the string topology of symmetric spaces of higher rank}, arXiv preprint arXiv:2212.09350 (2022).

\bibitem[Nae21]{naef:2021}
Florian Naef, \emph{The string coproduct" knows" reidemeister/whitehead torsion}, arXiv preprint arXiv:2106.11307 (2021).

\bibitem[NRW22]{naef:2022riveraWahl}
Florian Naef, Manuel Rivera, and Nathalie Wahl, \emph{String topology in three flavours}, arXiv preprint arXiv:2203.02429 (2022).

\bibitem[Oan15]{oancea:2015}
Alexandru Oancea, \emph{Morse theory, closed geodesics and the homology of free loop spaces}, Free loop spaces in geometry and topology, 2015, pp.~67--109.

\bibitem[Ste21]{stegemeyer:2021}
Maximilian Stegemeyer, \emph{On the string topology coproduct for lie groups}, arXiv:2109.10190, to appear in Homol. Homotopy Appl., 2021.

\bibitem[Zil77]{ziller:1977}
Wolfgang Ziller, \emph{The free loop space of globally symmetric spaces}, Inventiones Mathematicae \textbf{41} (1977), 1--22.

\end{thebibliography}
 \bibliographystyle{amsalpha}
 
\end{document}